\documentclass[12pt,a4paper]{article}
\usepackage{amssymb}
\usepackage{amsfonts}
\usepackage{amsmath}
\usepackage{amsthm}
\usepackage{enumerate}
\usepackage{srcltx}
\newtheorem{theorem}{Theorem}[section]
\newtheorem{proposition}{Proposition}[section]

\newtheorem{definition}{Definition}[section] 

\newtheorem{remark}{Remark}[section]
\newcommand{\wh}{\widehat}
\newcommand{\wt}{\widetilde}

\newcommand\cB{{\cal B}}

\newcommand\cX{{\cal X}}

\newcommand\cV{{\cal V}}

\newcommand\mes{\mbox{mes}}
\def\bbn{{\mathbb N}}

\def\bbr{{\mathbb R}}

\def\text#1{\hbox{#1}}
\def\proof{{\noindent \bf Proof. }}
\def\mes{{\bf mes}}
\def\endproof{\mbox{\ $\qed$}}

\def\E{{\bf E}}
\def\P{{\bf P}}
\def\C{{\bf C}}

\def\D{{\bf D}}
\def\B{{\bf B}}
\def\M{{\bf M}}
\def\A{{\bf A}}

\def\H{{\bf H}}

\def\x{{\bf x}}
\def\a{{\bf a}}
\def\b{{\bf b}}
\def\k{{\bf k}}
\def\g{{\bf g}}
\def\q{{\bf q}}
\def\Chi{{\bf 1}}
\def\d{\mathrm{d}}
\def\build #1_#2{\mathrel{\mathop{\kern 0pt #1}\limits_{#2}}}
\newcommand{\zs}[1]{{\mathchoice{#1}{#1}{\lower.25ex\hbox{$\scriptstyle#1$}}
{\lower0.25ex\hbox{$\scriptscriptstyle#1$}}}}

\numberwithin{equation}{section}
\def\proof{{\noindent \bf Proof. }}
\def\endproof{\mbox{\ $\qed$}}
\begin{document}
\title{
Geometric ergodicity for classes of homogeneous Markov chains
\thanks{The second author is partially supported by the RFFI-Grant  09-01-00172-a.}
}

\author{L. Galtchouk
\thanks{
IRMA,Department of Mathematics,
 Strasbourg University,
 7, rue R\'en\`e Descartes,
 67084, Strasbourg, France,
 e-mail: leonid.galtchouk@math.unistra.fr }
 \and
S. Pergamenshchikov\thanks{
 Laboratoire de Math\'ematiques Raphael Salem,
 Avenue de l'Universit\'e, BP. 12,
  Universit\'e de Rouen,
   F76801, Saint Etienne du Rouvray, Cedex France
and Department of Mathematics and Mechanics,Tomsk State University,
Lenin str. 36, 634041 Tomsk, Russia,
 e-mail: Serge.Pergamenchtchikov@univ-rouen.fr}
}

\maketitle

\begin{abstract}
The paper deals with non asymptotic computable bounds for the geometric convergence rate of homogeneous ergodic Markov processes.
Some sufficient conditions are stated for simultaneous geometric ergodicity
of Markov chain classes. This property is applied to 
 nonparametric estimating in ergodic diffusion processes.
\end{abstract}


\par


\bibliographystyle{plain}

\newpage

\section{Introduction}\label{sec:In}

In this paper we study a class of homogeneous Markov chains
$(\Phi)_\zs{n\ge 0}$ with values in some measurable
space $(\cX\,,\,\cB(\cX))$ defined  by a parametrized  family of transition probabilities
\begin{equation}\label{sec:In.1}
(\P^{\vartheta}_\zs{x})_\zs{x\in\cX\,,\,\vartheta\in \Theta}\,,
\end{equation}
where $\Theta$ is a parametric set for this family. For each
$\vartheta\in\Theta$ the sequence
$\Phi=(\Phi_\zs{n})_\zs{n\ge 0}$
is a homogeneous Markov chain defined on the measurable
space $(\cX\,,\,\cB(\cX))$ with a transition probability $\P^{\vartheta}$, i.e.
$$
\P(\Phi_\zs{1}\in\Gamma|\Phi_\zs{0}=x)=\P^{\vartheta}_\zs{x}(\Gamma)\,.
$$

Our main goal is to state
  geometric ergodicity  for this class simultaneously over all values of the parameter
$\vartheta\in\Theta$.

Geometric ergodicity is studied in a number of papers (see, for example, \cite{Ba}, \cite{MeTw}-\cite{NuTu}). We remind that a chain $(\Phi_\zs{n})_\zs{n\ge 0}$ on the space $(\cX,\cB(\cX))$ with a invariant measure $\pi$ is called
{\sl geometrically ergodic}
if there exist a $\cX\to [1,\infty[$ function $V(x)$ and some constants
$R>0$, $\kappa>0$ such that, for any $n\ge 1$,
\begin{equation}\label{sec:In.2}
\sup_\zs{x\in\cX}
\sup_\zs{0\le g\le V}
\frac{1}{V(x)}
|
\E_\zs{x}\,g(\Phi_\zs{n})
-
\pi(g)
|
\le R\, e^{-\kappa n}\,.
\end{equation}
  As we shall see later
(see, Definition~\ref{De.sec:Ug.1} below) the function $V$, providing the
 drift condition, is given by the Lyapunov functions (see, e.g. \cite{Kh} in the case of diffusion processes and \cite{FeTw}, \cite{KlPe1}, \cite{KlPe2} for Markov chains).
For this reason, in the sequel, we shall call such functions  by the
{\sl Lyapunov functions}.

The property  \eqref{sec:In.2} is useful
in applied problems related to the
identification of  stochastic systems, described by
 stochastic processes with dependent
values, in particular, governed by  stochastic difference or
stochastic differential equations.
Necessity of  simultaneous geometric ergodicity  appears in statistics,
when one studies a  minimax risk with respect to some family of distributions  related to
a statistical experiment.
 In particular, in this paper we shall apply simultaneous geometric
ergodicity  to nonparametric estimating the drift coefficient
in  the stochastic differential
equation (see \cite{GaPe1}):
\begin{equation}\label{sec:In.3}
\d y_\zs{t}=S(y_\zs{t})\,\d t +\sigma(y_\zs{t})\d W_\zs{t}\,,
\quad 0\le t\le T\,,
\end{equation}
where $S$ and $\sigma$ are unknown functions and $S$ has to be estimated from
observations $(y_\zs{t})_\zs{0\le t\le T}$. In studying  minimax risks for kernel estimators we need to use  geometric ergodicity for  the process \eqref{sec:In.3} simultaneously over all coefficients
$S$ and $\sigma$ from some functional class (see \eqref{sec:Mr.4} below).
In this case, $\vartheta=(S,\sigma)$ is the class parameter.

It is clear that to apply the property \eqref{sec:In.2} to some distribution family
we have to find some explicit expressions for the parameters $R>0$ and $\kappa>0$.
It is a well-known problem in the Markov Chain Monte Carlo (MCMC) theory, when a stopping rule for simulations is based on the accuracy of $n$-step approximations. Therefore, we need to find
 computable bounds in  \eqref{sec:In.2}. Note that some explicit expressions for $R$
and $\kappa$ were
calculated in \cite{Ba}, \cite{MeTw1}, \cite{RoRoSc}  and \cite{Ro}
for $\psi$-irreducible homogeneous Markov chains. These results are not applicable in our case because it is not clear what does it  mean $\psi$-irreducibility for a  class of parametrized Markov chains.
Note that
in \cite{Ba} and \cite{MeTw1} the parameters $R$ and $\kappa$ were obtained by making use of the Kendall theorem.
In \cite{Ro} these parameters are obtained through the
direct coupling method  for the Markov chain. To this end
the authors impose in \cite{RoRoSc}--\cite{Ro} some additional assumptions which are not satisfied in the case of the diffusion process
\eqref{sec:In.3} (see Remark~\ref{Re.sec:De.1} in  Section \ref{sec:De}).
Moreover, it should be noted that the upper bound in
\eqref{sec:In.2} given in \cite{RoRoSc} is calculated under
the assumption that the minorization condition holds on the whole state space.
This assumption is never true for the model \eqref{sec:In.3}.

In the paper we apply the coupling method to the renewal process generated by entrance
times of the process into the minorization set.
Note, that Meyn and Tweedie use the same approach in order 
to obtain  convergence results  (see, \cite{MeTw}, chapter 13). Their results imply the power convergence rate.
To obtain the geometric rate we make use of the Lyapunov functions method for the
related coupling process.

In order to explain the novelty of the method introduced in the paper, we give
the scheme of proving the property \eqref{sec:In.2}.
The first step consists in passing to a splitting
chain, which yields a chain with an accessible atom. Then, one makes use of the
{\sl Regenerative Decomposition} for splitting chains in order to evaluate the convergence rate (see \cite{Ki}, \cite{MeTw}).
Let us remind that the principal term in this decomposition gives the deviation
in the renewal theorem, which may be evaluated thanks to the
Kendall renewal theorem that provides a geometric convergence rate.
In our case  the same convergence rate is obtained
thanks to making use of the Lyapunov functions method for
the coupling renewal process (see Theorem~\ref{Th.sec:Cr.1} in Section~\ref{sec:Cr}).
 This upper bound enables us to find the explicit non asymptotic
  exponential upper bound in the ergodic theorem
for which we can find the supremum over the transition probability family in
\eqref{sec:In.1}. 

In this paper we find some   sufficient conditions
which provide simultaneous geometric ergodicity of the family \eqref{sec:In.1} over all values of the parameter $\vartheta$.
We check these conditions for the diffusion model \eqref{sec:In.3}. As corollary, we obtain explicit upper bounds for geometric
convergence rate in the ergodic theorem for diffusion processes. These bounds
may be used in the Monte Carlo technique to calculate some functionals of ergodic
diffusion processes. In that case one can replace these functionals
by the corresponding integrals with respect to the invariant density which has a simple explicit form. The accuracy of
this approximation is given by the explicit non asymptotic bounds in
the geometric convergence rate for diffusion processes.

The paper is organized as follows. In the next section the main results
are stated. Section~\ref{sec:De} provides the explicit formulas for parameters in the geometric convergence rate. Section~\ref{sec:Cr} is devoted to related coupling renewal processes.
In section~\ref{sec:Pr} geometric ergodicity is proved for a
parametrized class of homogeneous Markov chains. In section~\ref{sec:Ug}
we apply this property to stochastic differential equations.
 Some basic results on homogeneous Markov chains are given
in the Appendix.

\medskip

\section{Main results}\label{sec:Mr}

\noindent Assume now,  that the family of transition probabilities
$(\P^{\vartheta})_\zs{\vartheta\in\Theta}$ satisfies the
following conditions\\[1mm]

\medskip

\noindent $\H_\zs{1})$
{\em
There exist  $0<\delta<1$, 
 some set $C\in\cB_\zs(\cX)$ and some probability measure $\nu$ on $\cB(\cX)$
with $\nu(C)=1$ such that, for any $A\in\cB(\cX)$ with $\nu(A)>0$,
\begin{equation}\label{sec:Mr.2}
\inf_\zs{x \in C}\quad
\left(\inf_\zs{\vartheta \in\Theta}\,
\P^{\vartheta}(x,A)-\delta\nu(A)\right)\,>\,0\,.
\end{equation}
}
For the sequel we denote
\begin{equation}\label{sec:Mr.2-1}
\eta\,=\,\inf_\zs{x \in C}
\left(\inf_\zs{\vartheta \in\Theta}\,
\P^{\vartheta}(x,C)-\delta\right)\,.
\end{equation}
\noindent
\medskip

\medskip

\noindent $\H_\zs{2})$
{\em
There exist
$\cX\to [1,\infty)$ function $V$, some constants
$0<\rho<1$, $D\ge 1$, and a set $C$ from $\cB(\cX)$ such that
$$
V^{*}=\sup_\zs{x\in C}\,V(x)\,<\infty
$$
and, for any $x\in\cX$},
\begin{equation}\label{sec:Mr.1}
\sup_\zs{\vartheta\in\Theta}
\E^{\vartheta}_\zs{x}\left(V(\Phi_\zs{1})\right)\,
\le\,(1-\rho)V(x)\,+\,D\Chi_\zs{C}(x) \,.
\end{equation}

\medskip
\noindent Here $\E^{\vartheta}_\zs{x}$ means the expectation with respect to the transition probability $\P^{\vartheta}(x,\cdot)$.

\begin{remark}\label{Re.sec:Mr.1}
Condition $\H_\zs{2})$ is called the uniform
drift condition and that of $\H_\zs{1})$ is the uniform
minorization condition.
\end{remark}

\begin{theorem}\label{Th.sec:Mr.1}
Assume the conditions
$\H_\zs{1})$--$\H_\zs{2})$ hold true with the same set $C\in\cB(\cX)$.
Then, for each $\theta\in\Theta$, the chain $\Phi$ admits an invariant
distribution $\pi^{\vartheta}$ on $\cB(\cX)$. Moreover, for any $n\ge 2$,
\begin{equation}\label{sec:Mr.3}
\sup_\zs{\vartheta\in\Theta}\,
\sup_\zs{x\in\cX}\,
\sup_\zs{0< f\le V}
\frac{1}{V(x)}
\left|
\E^{\vartheta}_\zs{x}\,f(\Phi_\zs{n})
-
\int_\zs{\cX}\,f(z)\,\pi^{\vartheta}(\d z)
\right|
\le R^{*} \,e^{-\kappa^{*} n}\,,
\end{equation}
where the parameters $R^{*}=R^{*}(\rho,\delta,D,\eta,V^{*})$ and $\kappa^{*}=\kappa^{*}(\rho,\delta,D,\eta_\zs{1},V^{*})$ are given in
\eqref{sec:De.6}.
\end{theorem}

\medskip

\noindent Apply now to the process \eqref{sec:In.3}. To this end we have to introduce
 some functional class of functions $\vartheta=(S,\sigma)$.
First, for some $\x_\zs{*}\ge 1$, $M>0$ and $L>\beta>0$, we denote by $\cV_\zs{1}$
the class of functions $S$ from $\C^{1}(\bbr)$ such that
$$
\sup_\zs{|x|\le \x_\zs{*}}
\left(|S(x)|+|\dot{S}(x)|\right)\le M
$$
and
$$
-L\le \inf_\zs{|x|\ge \x_\zs{*}}\dot{S}(x)\le \sup_\zs{|x|\ge \x_\zs{*}}\dot{S}(x)
\le -\beta\,.
$$
Second, for some fixed reals
 $\sigma_\zs{min}>0$ and $\sigma_\zs{max}>\sigma_\zs{min}$, we denote by $\cV_\zs{2}$
the class of functions $\sigma$ from $\C^{2}(\bbr)$ such that,
for all $x\in\bbr$,
$$
\sigma_\zs{min}\le\,
\min\left(|\sigma(x)|\,,\,|\dot{\sigma}(x)|\,,\,|\ddot{\sigma}(x)|\right)
\le
\max\left(|\sigma(x)|\,,\,|\dot{\sigma}(x)|\,,\,|\ddot{\sigma}(x)|\right)
\le\sigma_\zs{max}\,.
$$
Finally, we set
\begin{equation}\label{sec:Mr.3-1}
\Theta=
\cV_\zs{1}\times \cV_\zs{2}
\,.
\end{equation}

\noindent
Note that (see, for example, \cite{GiSk}), for any function $\vartheta$ from
$\Theta$, the equation \eqref{sec:In.3} admits a unique strong solution, which is an ergodic process
having an invariant measure $\pi_\zs{\vartheta}$
with the invariant density $q_\zs{\vartheta}$ defined as
\begin{equation}\label{sec:Mr.3-2}
q_\zs{\vartheta}(x)=
\frac{\sigma^{-2}(x)\exp\{\int^{x}_\zs{0}\,S_\zs{1}(v)\d v\}}
{\int_{-\infty}^{+\infty}\,\sigma^{-2}(z)
\exp\{\int^{z}_\zs{0}\,S_\zs{1}(v)\d v\}
\d z}\,,
\end{equation}
where  $S_\zs{1}(v)= 2S(v)/\sigma^2(v)$.

\begin{theorem}\label{Th.sec:Mr.2}
For any $0<\epsilon\le 1/2$ and $t>0$,
\begin{equation}\label{sec:Mr.4}
\sup_\zs{\vartheta\in \Theta}
\sup_\zs{x\in\bbr}
\sup_\zs{0<g\le 1}\,
\frac{
\left|
\E^{\vartheta}_\zs{x}\,g(y_\zs{t})
-
\int_\zs{\bbr}g(x) q_\zs{\vartheta}(x)\d x
\right|}
{(1+x^2)^{\epsilon}}
\,\le R_\zs{\epsilon}\,e^{-\kappa_\zs{\epsilon} t}\,,
\end{equation}
where the parameters
$R_\zs{\epsilon}>0$ and
$\kappa_\zs{\epsilon}>0$
are given in \eqref{sec:De.10}.
\end{theorem}
\begin{remark}\label{Re.sec:Mr.2}
Note that the property \eqref{sec:Mr.4} is called simultaneous
geometric ergodicity. As is shown in Section~\ref{sec:Ug},
the function $(1+x^2)^{\epsilon}$
is a Lyapunov function. 
It should be said that in \cite{RoRo}, for the process
\eqref{sec:In.3} with a constant diffusion coefficient (i.e. for
$\sigma=1$), an exponential bound for deviation \eqref{sec:Mr.4} was obtained. 
The proof of that result was based on the coupling method applied directly to  
the diffusion process \eqref{sec:In.3} provided the existence of a Lyapunov function. In contrast with \cite{RoRo}, in our paper an explicit family of Lyapunov functions is given that is of help in applications.
\end{remark}

\begin{remark}\label{Re.sec:Mr.3}
It should be noted that the inequality \eqref{sec:Mr.4} may be applied
to Monte Carlo calculation of the expectation $\E^{\vartheta}_\zs{x}\,g(y_\zs{t})$. Indeed,
 the previous expectation can be replaced with the integral of $g$ with respect to the invariant density
\eqref{sec:Mr.3-2}. The precision of such approximation is given
in \eqref{sec:Mr.4}.
\end{remark}

\medskip
\medskip
\section{Computable bounds for geometric convergence rate}\label{sec:De}

In this section we introduce the parameters $R^{*}$ and $\kappa^{*}$
which make explicit the upper bound for geometric ergodicity in \eqref{sec:Mr.3}. 
For any $0<\gamma<1$, we denote

\begin{equation}\label{sec:De.1}
B^{*}=
\check{U}^{*}
\left(
1+
\frac{\check{U}^{*}V^{*}}
{1- (1-\delta \eta_\zs{1})^{\gamma}}
\right)
\,,
\end{equation}
where $\eta_\zs{1}=\eta/(1-\delta)$ and
$$
\check{U}^{*}=
\max\left(
\frac{1-\rho+D}{(1-\delta)(1-\rho)^{1-\gamma}\left(1-(1-\rho)^{\gamma}\right)}
\,,\,
\frac{V^{*}}{(1-\rho)^{1-\gamma}}
\right)
\,.
$$
We remind that the parameter $\eta$ appears in the condition $\H_\zs{1})$.
Moreover, we put
\begin{equation}\label{sec:De.1-1}
\left\{
\begin{array}{cl}
r_\zs{*}&
=\dfrac{(1-\gamma)^{2}|\ln(1-\rho)|\,|\ln(1-\delta\eta_\zs{1})|}
{\ln(\check{U}^{*}V^{*})+ |\ln(1-\delta\eta_\zs{1})|}
\,;
\\[8mm]
\wt{l}&
=
2
+
\left[
\dfrac{\ln(V^{*}B^{*})}{r_\zs{*}}
+
\dfrac{\ln\wt{q}(1-e^{-r_\zs{*}})^{-1}}{2r_\zs{*}}
\right]
\,,
\end{array}
\right.
\end{equation}
where
$$
\wt{q}=\frac{1-\check{B}^{*}_\zs{1}}{2}
\quad\mbox{and}\quad
\check{B}^{*}_\zs{1}
=
\min\left(e^{-r_\zs{*}}\,,\,
\frac{\delta}{e^{r_\zs{*}}-1+\delta\eta_\zs{1}}
\right)\,.
$$
Here $[a]$ means the integer part of a reel number $a>0$.
Next,
\begin{equation}\label{sec:De.3}
\wt{A}=
\dfrac{V^{*}B^{*}+1}{(1-e^{-r_\zs{*}})(1-e^{-\gamma r_\zs{*}})}
\quad\mbox{and}\quad
\wt{\varrho}_\zs{1}
=
\dfrac{(1-\gamma)^{2} \wt{r}\,|\ln(1-\wt{\epsilon})|}{\ln \wt{A}+r_\zs{*}\wt{l} 
+ |\ln(1-\wt{\epsilon})|}
\,,
\end{equation}
where 
$\wt{r}=|\ln(1-\wt{q})|$ and
$\wt{\epsilon}=\delta \eta_\zs{1}(1-\delta)^{\wt{l}-2}$.
Now we introduce the following parameter
\begin{equation}\label{sec:De.5}
\wt{A}_\zs{3}=
\frac{1-\gamma}{\gamma}
\left(
3 
+
2
\frac{
e^{r_\zs{*}}\,
\wt{A}\left(1+\wt{A}e^{r_\zs{*}\wt{l}}\right)}
{\left(1-(1-\wt{\epsilon})^{\gamma}\right)\left(e^{r_\zs{*}}-1\right)}
\right)\,
V^{*}B^{*}
+
e^{\wt{\kappa}}
\end{equation}
and
$$
\wt{\kappa}
=
\frac{(1-\gamma)\wt{\varrho}_\zs{1}r_\zs{*}}{\wt{\varrho}_\zs{1}+\ln V^{*}B^{*}}\,.
$$

\noindent We define the parameters $R^{*}$ and $\kappa^{*}$ in the Theorem~\ref{Th.sec:Mr.1} as

\begin{equation}\label{sec:De.6}
\left\{
\begin{array}{rl}
\kappa^{*}
=\kappa^{*}(\rho,\delta,D,V^{*})&=
\dfrac{(1-\gamma)\wt{\varrho}_\zs{1}r_\zs{*}}{2(\wt{\varrho}_\zs{1}+\ln V^{*}B^{*})}
\,;
\\[6mm]
R^{*}= R^{*}(\rho,\delta,D,V^{*})&=
2
\dfrac{(\wt{A}_\zs{2}+1)e^{\kappa^{*}}+1}{e^{\kappa^{*}}-1}\,V^{*}\,
\left( B^{*}\right)^{2}
\,.
\end{array}
\right.
\end{equation}
\noindent Further we define the upper bound in the ergodic Theorem~\ref{Th.sec:Mr.2}.
First of all, we take  the set $C$ 
 in the minorization condition  $\H_\zs{1})$ as the interval $C=[-K,K]$ 
 for some $K>0$ and  we chose  the  measure 
$\nu$ as the uniform distribution on $C$, i.e. for any measurable set $A$ from $\cB(\bbr)$,
\begin{equation}\label{sec:De.6-0}
\nu(A)=\,
\frac{1}{2 K}\,
\mes(A\cap [-K\,,\,K])\,,
\end{equation}
where $\mes(\cdot)$ is the Lebesgue measure on $\bbr$.

\noindent In order to define the threshold $\delta>0$ we need the quadratic function
\begin{equation}\label{sec:De.6-1}
\Omega_\zs{*}(z)=
\omega^{*}_\zs{1}\,z^{2}
+
\omega^{*}_\zs{2}
z
+
\omega^{*}_\zs{3}
\end{equation}
with $\omega^{*}_\zs{1}=4 (L\sigma_\zs{*})^{2}+ L\sigma_\zs{*}+1$,
$
\omega^{*}_\zs{2}=L\sigma_\zs{*}\sigma_\zs{max}
+
H^{*}_\zs{0}$ and
$\omega^{*}_\zs{3}=H^{*}_\zs{1}+2(H^{*}_\zs{0})^{2}$,
where
$$
H^{*}_\zs{0}=\frac{M+\sigma_\zs{*}}{\sigma_\zs{min}}\,,
\quad
H^{*}_\zs{1}=M(1+\sigma_\zs{*})+L+
2\sigma^{2}_\zs{*}
\quad\mbox{and}\quad
\sigma_\zs{*}=
\frac{\sigma_\zs{max}}{\sigma_\zs{min}}\,.
$$
We chose the parameter $\delta$ as
\begin{equation}\label{sec:De.6-2}
\delta_\zs{K}=
\frac{3 K e^{-\Omega_\zs{*}(\wt{K})}}{2\sqrt{2\pi} \sigma_\zs{max}}
\,,
\end{equation}
where  $\wt{K}=K/\sigma_\zs{min}$.
Now we set
\begin{equation}\label{sec:De.6-3}
K_\zs{0}=
\sqrt{\M_\zs{1}}
+8\sigma_\zs{min}+
8\check{\sigma}_\zs{*}
+
4\sqrt{\check{\sigma}_\zs{*}(\M_\zs{1}+\M_\zs{2})+16\check{\sigma}_\zs{*}^{2}}
\,,
\end{equation}
where $\check{\sigma}_\zs{*}
=\sigma_\zs{max}^{2}/(\beta(1-e^{-\beta}))$,
$$
\M_\zs{1}=
\frac{
(M+\beta\x_\zs{*})^2+\sigma^{2}_\zs{max}\beta}{\beta^{2}}
\quad\mbox{and}\quad 
\M_\zs{2}=
\M_\zs{1}
(1-e^{-\beta})
\,.
$$
We set 
\begin{equation}\label{sec:De.8}
\eta_\zs{K}=
1-\delta_\zs{K}-
\frac{4\check{\sigma}_\zs{*}
\left(K^{2}+\M_\zs{2}\right)}{\left(K^{2}-\M_\zs{1}\right)^{2}}
\,.
\end{equation}
Note, that  for any $K\ge K_\zs{0}$
$$
\delta_\zs{K}
\le
\frac{3 \wt{K} e^{-\wt{K}^{2}}}{2\sqrt{2\pi}}
\le 
\frac{3}{4\sqrt{e \pi}}
\quad\mbox{and}\quad
\eta_\zs{K}\ge 
\,
\frac{3(\sqrt{e \pi}-1)}{4\sqrt{e \pi}}
\,.
$$ 
Therefore, Propositions~\ref{Pr.subsec:AMc.7}--\ref{Pr.subsec:AMc.8} imply the 
condition $\H_\zs{1})$ with the parameters $\delta_\zs{K}$ and  $\eta_\zs{K}$
defined in
\eqref{sec:De.6-2} and \eqref{sec:De.8}
for any $K\ge K_\zs{0}$. To check the condition $\H_\zs{2})$
 we set
\begin{equation}\label{sec:De.8-1}
b^{*}_\zs{0}=
\beta (2+\x_\zs{*})
(1+\x_\zs{*})
+
(M+3) (3+\x_\zs{*})
\quad\mbox{and}\quad
b^{*}_\zs{1}=\frac{b^{*}_\zs{0}}{2\beta}
\,.
\end{equation}
Propositions~\ref{Pr.sec:Ug.2-1}--\ref{Pr.sec:Ug.3}
imply, that for any $0<\epsilon\le 1/2$,  the diffusion
process \eqref{sec:In.3} 
satisfies 
the drift condition $\H_\zs{2})$
with $V(x)=(1+x^{2})^{\epsilon}$, $C=[-K_\zs{\epsilon},K_\zs{\epsilon}]$, 
\begin{equation}\label{sec:De.9}
\rho_\zs{\epsilon}=(1-\check{\epsilon})\,(1-e^{-2\epsilon\beta})
\quad\mbox{and}\quad
D_\zs{\epsilon}=
V^{*}_\zs{\epsilon}
e^{-2\epsilon\beta}+ b^{*}_\zs{1}(1-e^{-2\epsilon\beta})\,,
\end{equation}
where
$$
K_\zs{\epsilon}=K_\zs{0}
+
\sqrt{\left(\frac{b^{*}_\zs{1}}{\check{\epsilon}}\right)^{1/\epsilon}-1}
\quad\mbox{and}\quad
V^{*}_\zs{\epsilon}=(1+K^{2}_\zs{\epsilon})^{\epsilon}
\,.
$$
Now we define
\begin{equation}\label{sec:De.10}
\kappa_\zs{\epsilon}=
\kappa^{*}
(\rho_\zs{\epsilon},\delta_\zs{\epsilon},D_\zs{\epsilon},\eta_\zs{\epsilon},V^{*}_\zs{\epsilon})
\quad\mbox{and}\quad
R_\zs{\epsilon}=
R^{*}
(\rho_\zs{\epsilon},\delta_\zs{\epsilon},D_\zs{\epsilon},\eta_\zs{\epsilon},V^{*}_\zs{\epsilon})
\,,
\end{equation}
where the functions $\kappa^{*}$ and $R^{*}$ are defined in
\eqref{sec:De.6}, $\delta_\zs{\epsilon}=\delta_\zs{K_\zs{\epsilon}}$,
$\eta_\zs{\epsilon}=\eta_\zs{K_\zs{\epsilon}}$
and $D_\zs{\epsilon}=D_\zs{K_\zs{\epsilon}}$.

\begin{remark}\label{Re.sec:De.1}
Note that we can not apply the bounds for geometric convergence rate from the paper
\cite{Ro} (Theorem 12) and \cite{RoRoSc} (Theorem 8) since there the bounds were obtained under the 
condition
\begin{equation}\label{sec:De.11}
(1+K^{2})^{\epsilon}>\frac{2D_\zs{K}}{\rho_\zs{\epsilon}}\,.
\end{equation}
This condition is not satisfied in 
our case for sufficiently small values of the parameter $\epsilon>0$. 
In fact, we apply the bounds \eqref{sec:Mr.4} in \cite{GaPe4}  to point-wise estimating the drift 
coefficient $S$ under observations of the process \eqref{sec:In.3} at discrete times
$(t_\zs{k})_\zs{k\ge 1}$. The question of interest is the behavior 
of kernel estimators that are nonlinear functionals of observations. 
In order to study the non asymptotic estimating precision  one needs 
of some concentration inequalities \cite{GaPe3} based
on the bounds \eqref{sec:Mr.4} with the parameter $\epsilon>0$ no matter how small.
\end{remark}

\noindent In order to illustrate the behavior of the geometric rate, we suppose that the
parameters satisfy the following conditions:
\begin{equation}\label{sec:De.12}
\left\{
\begin{array}{rl}
\lim_\zs{\delta\to 0}\eta=1\,,
\quad&
\lim_\zs{\delta\to 0}\,\dfrac{1-\eta}{\delta}
=+\infty\,,\\[4mm]
\lim_\zs{\rho\to 1} \dfrac{D}{V^{*}}=0
\quad&\mbox{and}\quad
\lim_\zs{\rho\to 1}
\dfrac{\ln V^{*}}{|\ln (1-\rho)|}=1\,.
\end{array}
\right.
\end{equation}
It should remark that these conditions hold true for the process \eqref{sec:In.3}
with the parametric set \eqref{sec:Mr.3-1} as $L>\beta\to\infty$ and
$\check{\epsilon}\to 0$ for some fixed $\epsilon>0$.

\noindent It should be noted, that the geometric rate in \eqref{sec:In.2} obtained in \cite{MeTw1} satisfies
$$
\kappa^{*}=O(\delta^{8})\quad\mbox{as}\quad \delta\to 0\,.
$$
\noindent In \cite{Ro} under the  condition  \eqref{sec:De.11}
this rate is ``best``, i.e.
$$
\kappa^{*}=O(\delta)\quad\mbox{as}\quad \delta\to 0\,.
$$
\noindent 
Under the conditions \eqref{sec:De.12} the coefficient $\kappa^{*}$ defined in
\eqref{sec:De.6}
$$
\kappa^{*}=(c_\zs{1}+o(1))\frac{\delta^{13/2+\mu_\zs{0}(\gamma)}}{|\ln \delta|^{2}}\quad\mbox{as}\quad 
\delta\to 0\,,\,\rho\to 1\,,
$$
where $c_\zs{1}>0$ and $\mu_\zs{0}(\gamma)\to 0$ as $\gamma\to 0$.

\begin{remark}\label{Re.sec:De.2}
Note that the condition $L,\beta\to\infty$ on the drift function of the process
\eqref{sec:In.3} concerns the behavior of the function only outside of the interval
$[-\x_\zs{*}\,,\,\x_\zs{*}]$, i.e. outside of the informative part of the function $S$. 
Remind (see, \cite{GaPe2}), that the class \eqref{sec:Mr.3-1} is used to bound the function $S$
on the interval $[-\x_\zs{*}\,,\,\x_\zs{*}]$ and outside of the interval the conditions are imposed to preserve 
 ergodicity. 
\end{remark}

\begin{remark}\label{Re.sec:De.3}
It should remark that, unfortunately, the bounds from the papers \cite{Ba}, \cite{MeTw1} are not applicable,
in general case, to classes of Markov chains. Indeed, irreducibility is one of conditions providing 
the geometrical rate in \cite{MeTw1}. Therefore, in the case of a parametric Markov chain class, irreducibility
measure should depend on a class parameter. It is not clear, what will going on this measure when one takes the
supremum over the class parameter and how one needs to change the irreducibility condition in order to obtain
uniform bounds over the class parameter $\vartheta \in\Theta$ by using the proof in \cite{MeTw1}.
\end{remark}

\section{Coupling Renewal Method}\label{sec:Cr}

In this Section we shall obtain a non asymptotic upper bound with explicit constants
in the renewal theorem by making use of the coupling method.  The notions
used here  can be found in \cite{Fe}, \cite{KlPe1}, \cite{Li}.

Let $(Y_\zs{j})_\zs{j\ge 0}$
and $(Y^{\prime}_\zs{j})_\zs{j\ge 0}$ be two independent sequences of random variables  taking values in $\bbn$. Assume that the initial random variables
 $Y_\zs{0}$ and $Y^{\prime}_\zs{0}$
have   distributions $\a=(\a(k))_\zs{k\ge 0}$ and $\b=(\b(k))_\zs{k\ge 0}$,
respectively, i.e. for any $k\ge 0$,
$$
\P(Y_\zs{0}=k)=\a(k)
\quad\mbox{and}\quad
\P(Y^{\prime}_\zs{0}=k)=\b(k)\,.
$$
The sequences $(Y_\zs{j})_\zs{j\ge 1}$
and $(Y^{\prime}_\zs{j})_\zs{j\ge 1}$ are supposed to be the i.i.d. sequences with the same distribution $p=(p(k))_\zs{k\ge 0}$, i.e. for any $k\ge 0$,
$$
\P(Y_\zs{1}=k)=\P(Y^{\prime}_\zs{1}=k)=p(k)\,.
$$
We assume also that  $p(0)=\P(Y_\zs{1}=0)=0$, i.e. the sequences
$(Y_\zs{j})_\zs{j\ge 1}$ and  $(Y^{\prime}_\zs{j})_\zs{j\ge 1}$
take  values in $\bbn^{*}=\bbn\setminus \{0\}$.
Moreover, we suppose that the distributions $a$, $b$ and $p$ satisfy the following condition

\medskip

\noindent $\C)$ {\em There exists a real number $r>0$ such that}
\begin{equation}\label{sec:Cr.2}
\ln\,\left(
\max\,(
\E\,e^{r Y_\zs{0}}\,,\,\E\,e^{r Y_\zs{1}})
\right)
\le
\upsilon_\zs{*}(r)
\quad\mbox{and}\quad
\ln\left(
\E\,e^{r Y'_\zs{0}}
\right)
\le \upsilon^{'}_\zs{*}(r)
\,.
\end{equation}

\medskip

\noindent For any $n\ge 0$, we define the following stopping times
$$
t_\zs{n}=\inf\{k\ge 0\,:\,\sum^{k}_\zs{i=0}\,Y_\zs{i}>n\}
\quad\mbox{and}\quad
t^{\prime}_\zs{n}=
\inf\{k\ge 0\,:\,\sum^{k}_\zs{i=0}\,Y^{\prime}_\zs{i}>n\}
\,.
$$
Further, we set
\begin{equation}\label{sec:Cr.3}
W_\zs{n}=\sum^{t_\zs{n}}_\zs{j=0}\,Y_\zs{j}-n
\quad\mbox{and}\quad
W^{\prime}_\zs{n}=\sum^{t^{\prime}_\zs{n}}_\zs{j=0}\,Y^{\prime}_\zs{j}-n\,.
\end{equation}
It is easy to see that the sequences
$(W_\zs{n})_\zs{n\ge 1}$ and $(W^{\prime}_\zs{n})_\zs{n\ge 1}$
are homogeneous Markov  chains taking values in
$\bbn^{*}$ such that, for any $n$, $k$ and $l$ from $\bbn^{*}$,
\begin{align}\nonumber
\P\left(W_\zs{n}=k|W_\zs{n-1}=l\right)
&=
\P\left(W^{\prime}_\zs{n}=k|W^{\prime}_\zs{n-1}=l\right)\\[2mm] \label{sec:Cr.3-1}
&=p(k)\Chi_\zs{\{l=1\}}+
\Chi_\zs{\{k=l-1\}}\Chi_\zs{\{l\ge 2\}}\,.
\end{align}

\noindent Firstly we study the entrance times
$(s_\zs{k})_\zs{k\ge 0}$
of the chain $(W_\zs{n})_\zs{n\ge 0}$ to the state $\{1\}$
which are defined as $s_\zs{-1}=0$
and for $k\ge 0$
\begin{equation}\label{sec:Cr.4}
s_\zs{k}=\inf\{l\ge s_\zs{k-1}+1\,:\,W_\zs{l}=1\}\,.
\end{equation}
\noindent One can check directly that
 the stopping times $s_\zs{k}, k\ge 1,$  can be represented as
\begin{equation}\label{sec:Cr.4-1}
s_\zs{k}=s_\zs{0}+
\sum^{k}_\zs{j=1} \varsigma_\zs{j}\,.
\end{equation}
One can check directly that in this case $(\varsigma_\zs{j})_\zs{j\ge 1}$ are i.i.d. 
random variables independent of
$s_\zs{0}$ and,
for any $l\in\bbn^{*}$,
\begin{equation}\label{sec:Cr.5}
\P(\varsigma_\zs{1}=l)=\P(s_\zs{0}=l|W_\zs{0}=1)=p(l)\,,
\end{equation}
i.e. the random variables $(\varsigma_\zs{j})_\zs{j\ge 1}$ have the same distribution as $Y_\zs{1}$.
Now we study the properties of the stopping time $s_\zs{0}$.
\begin{proposition}\label{Pr.sec:Cr.1}
Assume that the condition $\C)$ holds.
Then
$$
\E\,e^{r s_\zs{0}}\le 3\,e^{\upsilon_\zs{*}(r)}\,.
$$
\end{proposition}
\proof
First of all, note that, for $k\ge 2$ and $l\ge 0$,
$$
\P(s_\zs{0}=l|W_\zs{0}=k)=\Chi_\zs{\{l=k-1\}}\,.
$$
Moreover, taking into account that, for any $k\ge 1$, 
\begin{equation}\label{sec:Cr.5-1}
\P(W_\zs{0}=k)=\a(0)p(k)+\a(k)\,, 
\end{equation}
we obtain that
\begin{align*}
\E\,e^{r s_\zs{0}}
&=\P(W_\zs{0}=1)\,\E\,e^{r Y_\zs{1}}
+
\sum^{\infty}_\zs{k=2}\,
e^{r(k-1)}\,\P(W_\zs{0}=k)\\[2mm]
&\le 
\E\,e^{r Y_\zs{1}}
+
\a(0)e^{-r}\,\E\,e^{r Y_\zs{1}}
+
e^{-r}\,\E\,e^{r Y_\zs{0}}\,\le 3\,e^{\upsilon_\zs{*}(r)}\,.
\end{align*}
Hence Proposition~\ref{Pr.sec:Cr.1}.
\endproof

\noindent Now, we introduce the embedded Markov chain $(Z_\zs{k})_\zs{k\ge 0}$ by
\begin{equation}\label{sec:Cr.6}
Z_\zs{k}=W^{\prime}_\zs{s_\zs{k}}
\end{equation}
and the corresponding entrance time to the state $\{1\}$:
\begin{equation}\label{sec:Cr.7}
\varpi=\inf\{k\ge 1\,:\,Z_\zs{k}=1\}\,.
\end{equation}
In order to study the property of this stopping time, we need of the following
notations
\begin{equation}\label{sec:Cr.8}
l_\zs{*}=l_\zs{*}(r)=
2
+\left[
\frac{\ln\left(e^{2\upsilon_\zs{*}(r)}(1-e^{-r})^{-1}
\q(r)\right)}{2r}
\right]
\,,
\end{equation}
where 
$\q(r)=(1-e^{\upsilon_\zs{1}(r)})/2$ and the parameter 
$\upsilon_\zs{1}(r)<0$ will be specified below.
Moreover, for any $0<\gamma\,,\,\epsilon_\zs{*}<1$, we set 
\begin{equation}\label{sec:Cr.9}
A^{*}_\zs{1}(r)=A^{*}(r)
\,
\frac{1+A^{*}(r)e^{r l_\zs{*}}}{1-(1-\epsilon_\zs{*})^{\gamma}}
\,,
\end{equation}
where
$$
A^{*}(r)=
\frac{1+Q^{*}(r)e^{rl_\zs{*}+\upsilon_\zs{*}(r)}}
{(1-\q(r))^{1-\gamma}\,\left(1-(1-\q(r))^{\gamma}\right)}
\quad\mbox{and}\quad
Q^{*}(r)=\frac{e^{\upsilon_\zs{*}(r)}}{1-e^{-r}}
\,.
$$

\medskip
\begin{proposition}\label{Pr.sec:Cr.2}
Assume the condition $\C)$. Then, for  any $0<\gamma<1$, for
any  $\epsilon_\zs{*}>0$ and $\upsilon_\zs{1}(r)\le -r$  for which
\begin{equation}\label{sec:Cr.10}
0<
\epsilon_\zs{*}
\le
\min_\zs{1\le j\le l_\zs{*}-1}\,p(j)
\quad\mbox{and}\quad
\ln \E e^{-rY_\zs{1}}\le \upsilon_\zs{1}(r)
\,,
\end{equation}
one has
$$
\E\,e^{\varrho_\zs{*} \varpi}\le \,Q^{*}_\zs{1}(r)\,A^{*}_\zs{1}(r)\,,
$$
where
$Q^{*}_\zs{1}(r)=e^{\upsilon_\zs{*}(r)}
+
e^{\upsilon^{'}_\zs{*}(r)}
+
Q^{*}(r)$ and
\begin{equation}\label{sec:Cr.9-0}
\varrho_\zs{*}
=\varrho_\zs{*}(r)
=
\frac{(1-\gamma)^{2} |\ln(1-\q(r))|\,|\ln(1-\epsilon_\zs{*})|}{\ln A^{*}(r)+rl_\zs{*} 
+ |\ln(1-\epsilon_\zs{*})|}
\,.
\end{equation}
\end{proposition}
\proof
Firstly we note that the sequence \eqref{sec:Cr.6} is a homogeneous Markov chain with values in
$\bbn^{*}$ such that, for any $m$ and $l$ from $\bbn^{*}$ and for any $k\ge 1$,
\begin{equation}\label{sec:Cr.10-0}
\P(Z_\zs{k}=m|Z_\zs{k-1}=l)=
\P(W'_\zs{\varsigma_\zs{1}-l}=m|Y'_\zs{0}=0)\,.
\end{equation}
Note, that $W'_\zs{k}=-k$ for $k<0$
under the condition $Y'_\zs{0}=0$.
Therefore, for any positive function $V$,
\begin{equation}\label{sec:Cr.10-1}
\E\,[ V(Z_\zs{1})|Z_\zs{0}=l]=
\E\,V(l-\varsigma_\zs{1})\,\Chi_\zs{\{\varsigma_\zs{1}<l\}}
+
\E\,Q(\varsigma_\zs{1}-l)\Chi_\zs{\{\varsigma_\zs{1}\ge l\}}\,,
\end{equation}
where $Q(n)=\E\left(V(W_\zs{n})|Y_\zs{0}=0\right)$. Using the distribution
\eqref{sec:Cr.3-1} yields
$$
Q(n)
\le
\E\left(V(W_\zs{n-1}-1)|Y_\zs{0}=0\right)
+ \E V(Y_\zs{1}) \P\left(W_\zs{n-1}=1|Y_\zs{0}=0\right)\,.
$$
Choosing now $V(x)=e^{r x}$ one has
$$
Q(n)
\le
e^{-r} Q(n-1)
+ \E V(Y_\zs{1})\,.
$$
From the last inequality, taking into account that $Q(0)=\E V(Y_\zs{1})$,
 it follows that, for any $n\ge 1$,
\begin{equation}\label{sec:Cr.10-2}
Q(n)
\le
e^{-nr} Q(0)
+ \E V(Y_\zs{1})
\sum^{n}_\zs{j=1}\,e^{-(n-j)r}
\le 
Q^{*}(r)
\,,
\end{equation}
where the upper bound $Q^{*}(r)$ is defined in \eqref{sec:Cr.8}.
This implies that the last term in \eqref{sec:Cr.10-1} can be estimated as
$$
\E\,Q(\varsigma_\zs{1}-l)\Chi_\zs{\{\varsigma_\zs{1}\ge l\}}
\le Q^{*}(r) e^{\upsilon_\zs{*}(r)-rl}\,.
$$
Therefore,
$$
\frac{\E\,\left( V(Z_\zs{1})|Z_\zs{0}=l\right)}{V(l)}
\le 
e^{\upsilon_\zs{1}(r)}
+Q^{*}(r)e^{\upsilon_\zs{*}(r)-2rl}
\,.
$$
By making use of the definition of $l_\zs{*}$ in \eqref{sec:Cr.8}, we obtain
\begin{equation}\label{sec:Cr.11}
\sup_\zs{l\ge l_\zs{*}}
\frac{\E\left(V(Z_\zs{1})|Z_\zs{0}=l\right)}{V(l)}
\le
\frac{1+e^{\upsilon_\zs{1}(r)}}{2}
=
1-\q(r)<1\,.
\end{equation}
Moreover, for any $1\le l\le l_\zs{*}$,
$$
\E\left(V(Z_\zs{1})|Z_\zs{0}=l\right)
\le (1-\rho)V(l)+e^{\upsilon_\zs{1}(r)}V(l)
+Q^{*}(r)e^{\upsilon_\zs{*}(r)-rl}\,,
$$
i.e. the chain $(Z_\zs{k})_\zs{k\ge 1}$ satisfies the condition \eqref{subsec:AMc.2}
in the Appendix with
$$
C=\{1,\ldots,l_\zs{*}-1\}
\quad
\mbox{and}
\quad
D= e^{rl_\zs{*}+\upsilon_\zs{*}(r)}
Q^{*}(r)
\,.
$$
Therefore,  by Proposition~\ref{Pr.subsec:AMc.1}
for  $a_\zs{*}=-(1-\gamma)\ln (1-\q(r))$,
one gets
$$
\sup_\zs{l\ge 1}
\frac{U_\zs{C}(l,a_\zs{*},V)}{V(l)}
\le
\,A^{*}(r)
\,,
$$
where the upper bound $A^{*}(r)$ is given in
\eqref{sec:Cr.9}.
Moreover, from \eqref{sec:Cr.10} and \eqref{sec:Cr.10-0}
we get that, for  $2\le l\le l_\zs{*}-1$,
\begin{align*}
\P(Z_\zs{1}=1|Z_\zs{0}=l)
\ge p(l-1)\ge \epsilon_\zs{*}
\,.
\end{align*}
Therefore,
putting in Proposition~\ref{Pr.subsec:AMc.3} $\k_\zs{*}=\epsilon_\zs{*}$,
$a=\varrho_\zs{*}$ defined in \eqref{sec:Cr.9-0} and the set
$B=\{1\}$ we obtain that, for any $l\ge 1$,
$$
\E\,\left(e^{\varrho_\zs{*} \varpi}|Z_\zs{0}=l\right)
\le V(l) A^{*}_\zs{1}(r)\,,
$$
where the parameters $\varrho_\zs{*}$ and $A^{*}_\zs{1}(r)$
are defined in \eqref{sec:Cr.9} and \eqref{sec:Cr.9-0}. 
This upper bound implies 
$$
\E\,\left(e^{\varrho_\zs{*} \varpi}\right)
\le  A^{*}_\zs{1}(r)
\,\E\,V(Z_\zs{0}) 
\,.
$$
Moreover, note now that
$$
\E\,V(Z_\zs{0}) =\sum^{\infty}_\zs{j=1}\,
\E\,\left(V(Z_\zs{0})|W^{'}_\zs{0}=j\right)\,
\P(W^{'}_\zs{0}=j)\,.
$$
Similarly to \eqref{sec:Cr.10-1}
we obtain
$$
\E\,\left(V(Z_\zs{0})|W^{'}_\zs{0}=j\right)
=\E\,V(j-s_\zs{0})\,\Chi_\zs{\{s_\zs{0}<j\}}
+
\E\,Q(s_\zs{0}-j)\Chi_\zs{\{s_\zs{0}\ge j\}}\,.
$$
Using here the inequality \eqref{sec:Cr.10-2} yields
$$
\E\,\left(V(Z_\zs{0})|W^{'}_\zs{0}=j\right)\le V(j)+Q^{*}(r)\,.
$$
Moreover, similarly to \eqref{sec:Cr.5-1} we obtain that, for any $j\ge 1$,
$$
\P(W'_\zs{0}=j)=\b(0) p(j)+\b(j)\,.
$$
Thus,
$$
\E V(Z_\zs{0})\le \b(0)
\E V(Y_\zs{1})
+
\E V(Y^{'}_\zs{0})
+
Q^{*}(r)
$$
and we come to
the inequality
\eqref{sec:Cr.10}.
Hence Proposition~\ref{Pr.sec:Cr.2}.
\endproof

\begin{proposition}\label{Pr.sec:Cr.2-1}
 Assume  the condition $\C)$.
Then, for any $\epsilon_\zs{*}>0$ satisfying
the condition \eqref{Pr.sec:Cr.2} and for any $0<\gamma<1$,  
there exists $\kappa>0$ such that
\begin{equation}\label{sec:Cr.13}
\E e^{\kappa s_\zs{\varpi}}
\le
A^{*}_\zs{2}(r)
\,,
\end{equation}
where 
$$
A^{*}_\zs{2}(r)=
\frac{(1-\gamma)
\left(
3 e^{2\upsilon_\zs{*}(r)}
+Q^{*}_\zs{1}(r)
A^{*}_\zs{1}(r)
\right)
}{\gamma}
\quad\mbox{and}\quad
\kappa=\kappa(r)=
\frac{(1-\gamma)\varrho_\zs{*}r}{\varrho_\zs{*}+\upsilon_\zs{*}(r)}\,,
$$
the coefficients $A^{*}_\zs{1}(r)$ and $\varrho_\zs{*}$
are defined in 
\eqref{sec:Cr.9}
and
\eqref{sec:Cr.9-0}.
\end{proposition}
\proof Indeed, we have
\begin{align*}
 \E\,e^{\kappa s_\zs{\varpi}}&=
\kappa\,
\int^{\infty}_\zs{0}\,
e^{\kappa t}\,\P(s_\zs{\varpi}>t)\,\d t\\[2mm]
&\le 
\kappa\,
\int^{\infty}_\zs{0}\,
e^{\kappa t}\,
\left(
\P(s_\zs{N}>t)\,\d t
+
\P(\varpi>N)
\right)
\,\d t\,,
\end{align*}
where $N=N(t)=1+[\vartheta t]$, and $\vartheta$ is some positive parameter which will be 
chosen later. Note now that, for $0<\vartheta<r/\upsilon_\zs{*}(r) $,
$$
\P(s_\zs{N}>t)\le 3 e^{\upsilon_\zs{*}(r) (N+1)-rt}
\le 3 e^{2\upsilon_\zs{*}(r)}
\, e^{-(r-\upsilon_\zs{*}(r) \vartheta)t}\,.
$$
Moreover, due to Proposition~\ref{Pr.sec:Cr.2}
$$
\P(\varpi>N)\le Q^{*}_\zs{1}(r)\,A^{*}_\zs{1}(r)\,e^{-N\varrho_\zs{*}}
\le Q^{*}_\zs{1}(r)\,
A^{*}_\zs{1}(r)\,e^{-\vartheta \varrho_\zs{*}t}\,.
$$
Therefore, denoting
$$
\iota_\zs{*}(\vartheta)=\min
\left(
(r-\upsilon_\zs{*}(r) \vartheta)\,,\,
\vartheta \varrho_\zs{*}
\right)\,,
$$
one gets
$$
 \E\,e^{\kappa s_\zs{\varpi}}\le
\kappa\,\left(
3 e^{2\upsilon_\zs{*}(r)}
+Q^{*}_\zs{1}(r)\,
A^{*}_\zs{1}(r)
\right)\,
\int^{\infty}_\zs{0}\,e^{-(\iota_\zs{*}(\vartheta)-\kappa)t}\,
\d t\,.
$$
Maximizing now $\iota_\zs{*}(\vartheta)$ yields
$$
\max_\zs{0<\vartheta<\upsilon_\zs{*}(r)/r}\,\iota_\zs{*}(\vartheta)=
\iota_\zs{*}(\vartheta_\zs{max})
=\frac{r\varrho_\zs{*}}{\varrho_\zs{*}+\upsilon_\zs{*}(r)}
\,,
\quad
\vartheta_\zs{max}=\frac{r}{\varrho_\zs{*}+\upsilon_\zs{*}(r)}\,.
$$
Therefore, choosing now $\vartheta=\vartheta_\zs{max}$ and 
$\kappa=(1-\gamma)\iota_\zs{*}(\vartheta_\zs{max})$, we come to the inequality 
\eqref{sec:Cr.13}. Hence Proposition~\ref{Pr.sec:Cr.2-1}. \endproof

\noindent Let us define the renewal sequence $(u(n))_\zs{n\ge 0}$ as follows
\begin{equation}\label{sec:Cr.12}
u(n)=\sum^{\infty}_\zs{j=0}\,p^{*j}(n)
\,,
\end{equation}
where $p^{*j}$ denotes the $j$th convolution power.
For $j=0$ we set $p^{0}(n)=1$ for $n=0$ and $p^{0}(n)=0$ for $n\ge 1$.
We remind that, for two sequences $(\a(j))_\zs{j\ge 0}$ and $(u(j))_\zs{j\ge 0}$, the convolution
sequence $(\a*u(j))_\zs{j\ge 0}$
is
defined for any $j\ge 0$ as
$$
\a*u(j)=\sum^{j}_\zs{i=0}\,\a(i)u(j-i)\,.
$$

\bigskip

\begin{proposition}\label{Pr.sec:Cr.3}
 Assume that the condition $\C)$ holds and there exists $\epsilon_\zs{*}>0$ satisfying
the inequality \eqref{sec:Cr.10}.
Then, for any $0<\gamma<1$ and $n\ge 2$,
$$
\left|
\Delta(n)
\right|\,
\le A^{*}_\zs{2}\,e^{-\kappa n}\,,
$$
where $\Delta(n)=\a*u(n)
-
\b*u(n)$, the coefficients $\kappa$ and $A^{*}_\zs{2}$
are given in \eqref{sec:Cr.13}.
\end{proposition}
\proof
Obviously, that for $n\ge 1$,
$$
\a*u(n)=\P\left(\cup_\zs{j=0}
\left\{\sum^{j}_\zs{i=0}Y_\zs{i}=n\right\}\right)
=\P(W_\zs{n-1}=1)
$$
and
$$
\b*u(n)=
\P\left(\cup_\zs{j=0}\left\{\sum^{j}_\zs{i=0}\,Y^{\prime}_\zs{i}=n
\right\}\right)
=\P(W^{\prime}_\zs{n-1}=1)\,.
$$
Therefore,
$$
\Delta(n)=
\P(W_\zs{n-1}=1\,,\,W^{\prime}_\zs{n-1}\ge 2)
-
\P(W^{\prime}_\zs{n-1}=1\,,\,W_\zs{n-1}\ge 2)\,.
$$
Now, we introduce the ``coupling'' stopping time $\tau$
as
$$
\tau=\inf\{k\ge 1\,:\,(W_\zs{k},W^{\prime}_\zs{k})=(1,1)\}\,.
$$
Note that, for any $n\ge 2$, by the Markov property for the chain
$(W_\zs{k},W^{\prime}_\zs{k})_\zs{k\ge 1}$, one has
\begin{align*}
\P(W_\zs{n}=1\,,\,W^{\prime}_\zs{n}\ge 2,\tau\le n-1)
&=\sum^{n-1}_\zs{k=1}\,
\P(W_\zs{n}=1\,,\,W^{\prime}_\zs{n}\ge 2,\tau=k)\\[2mm]
&=
\sum^{n-1}_\zs{k=1}\,\P(\tau=k)\,v_\zs{n-k}\,.
\end{align*}
where $v_\zs{k}=\P(W_\zs{k}=1\,|\,W_\zs{0}=1)\P(W_\zs{k}\ge 2\,|\,W_\zs{0}=1)$.
Similarly, one gets
$$
\P(W^{\prime}_\zs{n}=1\,,\,W_\zs{n}\ge 2,\tau\le n-1)
=\sum^{n-1}_\zs{k=1}\,\P(\tau=k) v_\zs{n-k}\,.
$$
This implies that
$$
\Delta(n)=
\alpha_\zs{1}(n-1)-\alpha_\zs{2}(n-1)\,,
$$
where $\alpha_\zs{1}(n)=
\P(W_\zs{n}=1\,,\,W^{\prime}_\zs{n}\ge 2\,,\tau> n)$
and
$$
\alpha_\zs{2}(n)=\P(W^{\prime}_\zs{n}=1\,,\,W_\zs{n}\ge 2\,,\tau> n)\,.
$$
Therefore, for any $n\ge 2$,
\begin{align*}
|\Delta(n)|\le\,
\max\left(\alpha_\zs{1}(n-1)\,,\,\alpha_\zs{2}(n-1)\right)
\le
\P(\tau> n)\,.
\end{align*}
Taking into account that $\tau\le s_\zs{\varpi}$
a.s., we obtain
$$
|\Delta(n)|\le \P(s_\zs{\varpi}>n)\le
\,e^{-\kappa n}\,
\E e^{\kappa s_\zs{\varpi}}\,.
$$
Proposition~\ref{Pr.sec:Cr.2-1} implies
the upper bound \eqref{sec:Cr.13}. Hence Proposition~\ref{Pr.sec:Cr.3}
\endproof

\medskip

\begin{theorem}\label{Th.sec:Cr.1}
Assume that
there exists $r>0$ such that
$$
\ln\,\E\,e^{rY_\zs{1}}\le \upsilon_\zs{*}(r)\,.
$$
Then, for any $0<\gamma<1$, $n\ge 2$
and $\epsilon_\zs{*}>0$ satisfying
the inequality \eqref{sec:Cr.10},
$$
\left|
u(n)
-
\frac{1}{\E\, Y_\zs{1}}
\right|\,
\le 
A^{*}_\zs{3}(r)
\,e^{-\kappa n}\,,
$$
where
\begin{equation}\label{sec:Cr.15}
A^{*}_\zs{3}(r)
=
\frac{1-\gamma}{\gamma}
\left(
3
+2
\frac{e^{r}\,A^{*}(r)\,
\left(1+A^{*}(r)e^{rl_\zs{*}}\right)
}{\left(1-(1-\epsilon_\zs{*})^{\gamma}\right)(e^{r}-1)}
\right)\,e^{\upsilon_\zs{*}(r)}
\end{equation}
and the parameter $\kappa>0$ is defined in \eqref{sec:Cr.13}.
\end{theorem}
\proof
We obtain the inequality \eqref{sec:Cr.15} through Proposition~\ref{Pr.sec:Cr.3} in which
we choose $\a(0)=1$ with $\a(j)=0$ for $j\ge 1$. Moreover, we choose
the distribution $(\b(j))_\zs{j\ge 0}$  as
$$
\b(j)=\frac{1}{\E Y_\zs{1}}\,\P(Y_\zs{1}>j)
=\frac{1}{\E Y_\zs{1}}\,\sum^{\infty}_\zs{i=j+1}\,p(i)\,.
$$
It is easy to see directly that, for any $j\ge 1$,
$$
\b*u(j)=\frac{1}{\E Y_\zs{1}}\,.
$$
Note now, that through  the condition of this theorem we obtain 
$$
e^{r n}\,\b(n)=
\frac{e^{r n}\P(Y_\zs{1}>n)}{\E Y_\zs{1}}
\le \frac{\E e^{r Y_\zs{1}}\Chi(Y_\zs{1}>n)}{\E Y_\zs{1}}
\to
0\,,\quad\mbox{as}\quad n\to\infty\,.
$$
Therefore, the summing by parts yields
$$
\sum_\zs{j\ge 0}\,e^{rj}\b(j)=
\frac{\E\,e^{rY_\zs{1}}-1}{(e^{r}-1)\E Y_\zs{1}}
\le \frac{e^{\upsilon_\zs{*}(r)}-1}{e^{r}-1}
:=e^{\upsilon^{'}_\zs{*}(r)}\,.
$$
and Proposition~\ref{Pr.sec:Cr.3} implies the
inequality \eqref{sec:Cr.15}.
Hence Theorem~\ref{Th.sec:Cr.1}.
\endproof

\medskip

\section{Proof of Theorem~\ref{Th.sec:Mr.1}}\label{sec:Pr}

First we fix some $0<\gamma<1$ and we set $a_\zs{1}=-(1-\gamma)\ln (1-\rho)$.
We start with studying the properties of the function
$$
U^{\vartheta}_\zs{C}(x,a_\zs{1},V)=
\E^{\vartheta}_\zs{x}\,\sum^{\tau_\zs{C}}_\zs{j=1}\,
e^{a_\zs{1} j}\,V(\Phi_\zs{j})
$$
where $\tau_\zs{C}=\inf\{n\ge 1\,:\,\Phi_\zs{n}\in C\}$.
The condition $\H_\zs{2})$
and Proposition~\ref{Pr.subsec:AMc.1} imply immediately
that, for any $0<\gamma<1$,
\begin{equation}\label{sec:Pr.1}
\sup_\zs{\vartheta\in\Theta}\,\sup_\zs{x\in\cX}\,
\frac{U^{\vartheta}_\zs{C}(x,a_\zs{1},V)}{V(x)}
\le \frac{1-\rho+D}{(1-\rho)^{1-\gamma}(1-(1-\rho)^{\gamma})}:=U^{*}
\,.
\end{equation}
Now, we introduce a splitting chain family as in
  \cite{MeTw}, p. 108 (see also \cite{Nu}).
 We set
$\check{\cX}=\cX\times \{0,1\}$, $\cX_\zs{0}=\cX\times\{0\}$ and
$\cX_\zs{1}=\cX\times\{1\}$. Let $\cB(\cX_\zs{i})$ be the $\sigma-$fields generated
by the set $A_\zs{i}=A\times\{i\}$ with $A\in\cB(\cX)$, $i=0,1$.
In the sequel we will denote by
$<\check{x}>_\zs{i}$ the $i$th component of $\check{x}\in\check{\cX}$.
It is clear, that $<\check{x}>_\zs{0}\in\cX$ and
$<\check{x}>_\zs{1}\in\{0,1\}$. Furthermore, we define the $\sigma$ - field
$\cB(\check{\cX})$ as a $\sigma-$field generated by
$\cB(\cX_\zs{0})\cup\cB(\cX_\zs{1})$ and for any measure
$\lambda$ on $\cB(\cX)$ we relate the measure $\lambda^{*}$ on $\cB(\check{\cX})$
as
$$
\lambda^{*}(A_\zs{0})=(1-\delta)\lambda(A\cap C)+\lambda(A\cap C^{c})
\quad\mbox{and}\quad
\lambda^{*}(A_\zs{1})=\delta\lambda(A\cap C)\,.
$$
Now, for each $\vartheta\in\Theta$, we introduce a homogeneous Markov chain
$(\check{\Phi}_\zs{n})_\zs{n\ge 0}$
by the following transition probabilities
\begin{equation}\label{sec:Pr.2}
\check{\P}^{\vartheta}(\check{x},\cdot)
=\left\{
\begin{array}{ccc}
\P^{\vartheta}(x,\cdot)^{*}\,,&\quad\mbox{if}\quad &\check{x}\in \cX_\zs{0}\setminus C_\zs{0}\,;\\[3mm]
\dfrac{
\P^{\vartheta}(x,\cdot)^{*}
-
\delta \nu^{*}(\cdot)
}{1-\delta}
\,,&
\quad\mbox{if}\quad &\check{x}\in C_\zs{0}\,;\\[5mm]
\nu^{*}(\cdot)\,,&
\quad\mbox{if}\quad &\check{x}\in \cX_\zs{1}\,.
\end{array}
\right.
\end{equation}
Note, that for any
$\check{x}\in\cX_\zs{1}$,
\begin{equation}\label{sec:Pr.2-1}
\check{\P}^{\vartheta}(\check{x},C_\zs{0}\cup C_\zs{1})
=
\nu^{*}(C_\zs{0}\cup C_\zs{1})=
\nu(C)=1\,.
\end{equation}
Obviously, that the set
$\alpha=C_\zs{1}$ is an accessible  atom for the chain
$(\check{\Phi}_\zs{n})_\zs{n\ge 1}$, i.e. for any positive
$\check{\cX}\to\bbr$
function
$g$
$$
\check{\E}^{\vartheta}_\zs{\check{x}}\,g(\check{\Phi}_\zs{1})=
\check{\E}^{\vartheta}_\zs{\check{y}}\,g(\check{\Phi}_\zs{1})\,,
\quad\mbox{for any}
\quad \check{x}\,,\,\check{y}\in \alpha\,.
$$

This implies directly that, for any nonnegative  random variable $\xi$ measurable with respect to the
$\sigma-$field generated by the chain
$(\check{\Phi}_\zs{n})_\zs{n\ge 1}$, one has
$$
\check{\E}^{\vartheta}_\zs{\check{x}}\,\xi=
\check{\E}^{\vartheta}_\zs{\check{y}}\,\xi
\quad\mbox{for any}
\quad \check{x}\,,\,\check{y}\in \alpha\,.
$$
In the sequel we denote by $\check{\E}^{\vartheta}_\zs{\alpha}(\cdot)$ the such expectations.
Moreover, one can check directly that the chain $(\check{\Phi}_\zs{n})_\zs{n\ge 1}$ is
$\nu^{*}$-irreducible. 
 Next, for any set $\check{C}$ from $\cB(\check{\cX})$ we introduce
the corresponding entrance time
\begin{equation}\label{sec:Pr.3}
\check{\tau}_\zs{\check{C}}=\inf\left\{n\ge 1\,:\,\check{\Phi}_\zs{n}\in
\check{C}
\right\}
\end{equation}
and the corresponding entrance function
\begin{equation}\label{sec:Pr.5}
\check{U}^{\vartheta}_\zs{\check{C}}(\check{x},a_\zs{1},\check{V})=
\check{\E}^{\vartheta}_\zs{\check{x}}\,
\sum^{\check{\tau}_\zs{\check{C}}}_\zs{j=1}\,
e^{a_\zs{1} j}\,\check{V}(\check{\Phi}_\zs{j})\,,
\end{equation}
where $\check{V}(\check{x})=V(<\check{x}>_\zs{1})$.
By Proposition~\ref{Pr.subsec:AMc.5} we obtain that, for any
$x\in\cX$,
\begin{align*}
U^{\vartheta}_\zs{C}(x,a_\zs{1},V)&=(1-\delta)
\check{U}^{\vartheta}_\zs{C_\zs{0}\cup C_\zs{1}}(x_\zs{0},a_\zs{1},\check{V})\,
\Chi_\zs{\{x\in C\}}
+
\check{U}^{\vartheta}_\zs{C_\zs{0}\cup C_\zs{1}}(x_\zs{0},a_\zs{1},\check{V})\,
\Chi_\zs{\{x\in C^{c}\}}\\[2mm]
&+\delta
\check{U}^{\vartheta}_\zs{C_\zs{0}\cup C_\zs{1}}(x_\zs{1},a_\zs{1},\check{V})\,
\Chi_\zs{\{x\in C\}}\,,
\end{align*}
where
$x_\zs{i}=(x,i)$ for $i=0,1$.
Note now, that due to the property \eqref{sec:Pr.2-1},
$$
\check{U}^{\vartheta}_\zs{C_\zs{0}\cup C_\zs{1}}(x_\zs{1},a_\zs{1},\check{V})
=\check{\E}^{\vartheta}_\zs{x_\zs{1}}\,
e^{a_\zs{1}}\,\check{V}(\check{\Phi}_\zs{1})
\le e^{a_\zs{1}}\,V^{*}
=\frac{V^{*}}{(1-\rho)^{1-\gamma}}\,.
$$
Therefore, using the upper bound \eqref{sec:Pr.1}
and the coefficient $\check{U}^{*}$
given in \eqref{sec:De.1}, we obtain
$$
\sup_\zs{\vartheta\in\Theta}\,\sup_\zs{\check{x}\in\check{\cX}}\,
\frac{\check{U}^{\vartheta}_\zs{C_\zs{0}\cup C_\zs{1}}(\check{x},a_\zs{1},\check{V})}{\check{V}(\check{x})}
\le \check{U}^{*}
\,.
$$
Note now that, for $\check{x}\in C_\zs{0}$
by the definition \eqref{sec:Pr.2}
$$
\check{\P}^{\vartheta}(\check{x},\alpha)=\check{\P}^{\vartheta}(\check{x},C_\zs{1})
=\delta\frac{\P^{\vartheta}(x,C)-\delta}{1-\delta}
\ge \delta \eta_\zs{1}\,,
$$
where the parameter $0<\eta_\zs{1}\le 1$
is defined in \eqref{sec:De.1}.
By making use of Proposition~\ref{Pr.subsec:AMc.3}
with $a_\zs{*}=a_\zs{1}=-(1-\gamma)\ln (1-\rho)$ and 
$\k_\zs{*}=\delta\eta_\zs{1}$, one gets
\begin{equation}\label{sec:Pr.6}
\sup_\zs{\vartheta\in\Theta}\,\sup_\zs{\check{x}\in\check{\cX}}\,
\frac{\check{U}^{\vartheta}_\zs{\alpha}(\check{x},r_\zs{*},\check{V})}{\check{V}(\check{x})}
\le B^{*}
\,,
\end{equation}
where
$B^{*}$ and $r_\zs{*}$
are given in \eqref{sec:De.1}
and \eqref{sec:De.1-1}.
Therefore, by Proposition~\ref{Pr.subsec:AMc.0},
 the chain $(\check{\Phi}_\zs{n})_\zs{n\ge 0}$
is ergodic for each $\vartheta\in \Theta$ with the invariant measure given as

\begin{equation}\label{sec:Pr.7}
\check{\pi}^{\vartheta}(\check{\Gamma})
=
\frac{1}{\check{\E}^{\vartheta}_\zs{\alpha} \check{\tau}_\zs{\alpha}}
\,\check{\E}^{\vartheta}_\zs{\alpha}\,
\sum^{\check{\tau}_\zs{\alpha}}_\zs{j=1}\,
\Chi_\zs{\{\check{\Phi}_\zs{j}\in\check{\Gamma}\}}\,.
\end{equation}
Now, for any $n\ge 2$, we define
$
\check{\iota}=\max\{1\le j\le n-1\,:\,\check{\Phi}_\zs{j}\in\alpha \}
$
and we put $\check{\iota}=0$ if $\check{\tau}_\zs{\alpha}\ge n$. Moreover, note that,
for any $\check{\cX}\to\bbr$
 function $f$ and any $n\ge 2$,
\begin{align*}
\check{\E}^{\vartheta}_\zs{\check{x}}\,f(\check{\Phi}_\zs{n})\,
\Chi_\zs{\{\check{\tau}_\zs{\alpha}<n\}}
&=\sum^{n-1}_\zs{j=1}\,\check{\E}^{\vartheta}_\zs{\check{x}}\,
f(\check{\Phi}_\zs{n})\,
\Chi_\zs{\{\check{\tau}_\zs{\alpha}\le j\}}
\Chi_\zs{\{\check{\iota}=j\}}
\\[2mm]
&=\sum^{n-1}_\zs{j=1}\,\check{\E}^{\vartheta}_\zs{\check{x}}\,
\Chi_\zs{\{\check{\tau}_\zs{\alpha}\le j\}}\,
\check{\E}^{\vartheta}_\zs{\check{x}}\,
\left(f(\check{\Phi}_\zs{n})\,
\Chi_\zs{\{\check{\iota}=j\}}
| \check{\Phi}_\zs{1},\ldots,\check{\Phi}_\zs{j}\right)\,.
\end{align*}
Now, taking into account that
$(\check{\Phi}_\zs{n})_\zs{n\ge 1}$
is a homogeneous Markov chain, we can calculate  the last conditional expectation
as follows
\begin{align}\nonumber
\check{\E}^{\vartheta}_\zs{\check{x}}\,
\left(f(\check{\Phi}_\zs{n})\,
\Chi_\zs{\{\check{\iota}=j\}}
| \check{\Phi}_\zs{1},\ldots,\check{\Phi}_\zs{j}\right)
&=
\Chi_\zs{\{\check{\Phi}_\zs{j} \in\alpha \}}
\check{\E}^{\vartheta}_\zs{\alpha}\,
\left(f(\check{\Phi}_\zs{n-j})\,
\Chi_\zs{
\{
\check{\Phi}_\zs{1} \notin\alpha\,
,\ldots,
\check{\Phi}_\zs{n-j-1}\notin \alpha
\}}\right)
\\[4mm]\label{sec:Pr.8}
&=
\Chi_\zs{\{\check{\Phi}_\zs{j} \in\alpha \}}
g_\zs{f,\alpha}(n-j)\,,
\end{align}
where 
$g_\zs{f,\alpha}(k)=\check{\E}^{\vartheta}_\zs{\alpha}\,f(\check{\Phi}_\zs{k})\,
\Chi_\zs{\{\check{\tau}_\zs{\alpha}\ge k\}}$.
By convention, we set $g_\zs{f,\alpha}(0)=0$. Therefore,
$$
\check{\E}^{\vartheta}_\zs{\check{x}}\,f(\check{\Phi}_\zs{n})\,
\Chi_\zs{\{\check{\tau}_\zs{\alpha}<n\}}=
\sum^{n}_\zs{j=1}\,
\check{\P}^{\vartheta}_\zs{\check{x}}\left(\check{\tau}_\zs{\alpha}\le j \right)\,
g_\zs{f,\alpha}(n-j)
=h_\zs{\check{x}}*g_\zs{f,\alpha}(n)\,,
$$
where $h_\zs{\check{x}}(0)=0$ and, for $j\ge 1$,
$$
h_\zs{\check{x}}(j)=
\check{\P}^{\vartheta}_\zs{\check{x}}\left(\check{\tau}_\zs{\alpha}\le j \right)
=\check{\P}^{\vartheta}_\zs{\check{x}}\left(
\check{\Phi}_\zs{j}
\in
\alpha
 \right)\,.
$$
Moreover, for $j\ge 1$,
\begin{align*}
h_\zs{\check{x}}(j)
=
\sum^{j}_\zs{l=1}\,
\check{\P}^{\vartheta}_\zs{\check{x}}\left(
\check{\tau}_\zs{\alpha}=l\,,\,
\check{\Phi}_\zs{j}
\in
\alpha
 \right)=
\sum^{j}_\zs{l=1}\,\a_\zs{\check{x}}(l)\,u(l-j)
=\a_\zs{\check{x}}*u(j)
\,,
\end{align*}
where
\begin{equation}\label{sec:Pr.9}
\a_\zs{\check{x}}(l)=\check{\P}^{\vartheta}_\zs{\check{x}}\left(
\check{\tau}_\zs{\alpha}=l
 \right)
\quad\mbox{and}\quad
u(l)=\check{\P}^{\vartheta}_\zs{\alpha}\left(
\check{\Phi}_\zs{l}
\in
\alpha
 \right)
\,.
\end{equation}
It is clear that $\a_\zs{\check{x}}(0)=0$ and $u(0)=1$, i.e.
$\a_\zs{\check{x}}*u(0)=0$. This implies that
$$
h_\zs{\check{x}}(j)=\a_\zs{\check{x}}*u(j)\,,
$$
for all $j\ge 0$.
Finally, taking into account that $\a_\zs{\check{x}}*u*g_\zs{f,\alpha}(1)=0$, we obtain that for any $n\ge 1$,
\begin{equation}\label{sec:Pr.10}
\check{\E}^{\vartheta}_\zs{\check{x}}\,f(\check{\Phi}_\zs{n})\,
\Chi_\zs{\{\check{\tau}_\zs{\alpha}<n\}}=
\a_\zs{\check{x}}*u*g_\zs{f,\alpha}(n)\,.
\end{equation}
Note that the sequence $(u(n))_\zs{n\ge 0}$ is a renewal sequence, i.e.
$$
u(n)=\sum^{\infty}_\zs{j=0}\,p^{*j}(n)
\quad\mbox{with}\quad
p(k)=
\check{\P}^{\vartheta}_\zs{\alpha}(\check{\tau}_\zs{\alpha}=k)\,.
$$
Now, we set $Y_\zs{0}=\check{\tau}_\zs{\alpha}$ and
$
Y_\zs{j}=\inf\{j\ge Y_\zs{j-1}+1\,:\,\check{\Phi}_\zs{j}\in
\alpha\}$.
One can check directly that 
$(Y_\zs{j})_\zs{j\ge 1}$ is i.i.d. sequence with the distribution
$(p(k))_\zs{k\ge 1}$, i.e. $u(n)_\zs{n\ge 0}$ is the renewal function
 for the sequence $(Y_\zs{j})_\zs{j\ge 1}$.
We denote
$$
\omega(n)=\left|
u(n)-\frac{1}{\check{\E}^{\vartheta}_\zs{\alpha}\,Y_\zs{1}}
\right|\,.
$$
We estimate this term by Theorem~\ref{Th.sec:Cr.1}. First we have to check
the condition
$\C_\zs{1})$ uniformly over the parameter $\vartheta\in\Theta$, i.e.
to show that, for any $j\ge 1$,
\begin{equation}\label{sec:Pr.11}
\inf_\zs{\vartheta\in\Theta}\,
\check{\P}^{\vartheta}_\zs{\alpha}(\check{\tau}_\zs{\alpha}=j)
\ge 
\delta\,\eta_\zs{1}(1-\delta)^{j-1}\,.
\end{equation}
Let us check this property for $j=1$. We remind that, by the condition
$\H_\zs{1})$, one has $\nu(C)=1$. Thus, the definition \eqref{sec:Pr.2}
implies
$$
\check{\P}^{\vartheta}_\zs{\alpha}(\check{\tau}_\zs{\alpha}=1)=
\nu^{*}(C_\zs{1})=\delta\ge \delta\,\eta_\zs{1}\,.
$$
To show the property \eqref{sec:Pr.11} for $j\ge 2$
note, that
$
\check{\P}^{\vartheta}_\zs{\alpha}(\cX_\zs{0})=\nu^{*}(\cX_\zs{0})=
1- \delta$. Moreover, taking into account that
$\P^{\vartheta}(z,\cX_\zs{0})^{*}\ge 1-\delta$ we obtain
the same lower bound for the splitting distribution \eqref{sec:Pr.2}, i.e.
  for any $\check{z}\in \cX_\zs{0}$
\begin{align*}
\check{\P}^{\vartheta}(\check{z},\cX_\zs{0})&=
\Chi_\zs{\{\check{z}\in C_\zs{0}\}}
\left(
\frac{\P^{\vartheta}(z,\cX_\zs{0})^{*}- \delta\nu^{*}(\cX_\zs{0})}{1-\delta}
\right)
\,\\[2mm]
&
+
\Chi_\zs{\{\check{z}\in \cX_\zs{0}\setminus C_\zs{0}\}}
\P^{\vartheta}(z,\cX_\zs{0})^{*}
\ge 
(1-\delta)
\,.
\end{align*}
Similarly, for any $\check{z}\in \cX_\zs{0}$ we obtain
\begin{align*}
\check{\P}^{\vartheta}(\check{z},C_\zs{1})&=
\Chi_\zs{\{\check{z}\in C_\zs{0}\}}
\left(
\frac{\P^{\vartheta}(z,C_\zs{1})^{*}- \delta\nu^{*}(C_\zs{1})}{1-\delta}
\right)
\,\\[2mm]
&
+
\Chi_\zs{\{\check{z}\in \cX_\zs{0}\setminus C_\zs{0}\}}
\P^{\vartheta}(z,C_\zs{1})^{*}
\ge 
\delta
\frac{\P^{\vartheta}(z,C)-\delta}{1-\delta}
\ge \delta\eta_\zs{1}
\,.
\end{align*}
Now through the induction
we can show, that for any $j\ge 1$
$$
\check{\P}^{\vartheta}_\zs{\alpha}
\left(
\check{\Phi}_\zs{1}\in \cX_\zs{0}\,,\ldots,
\check{\Phi}_\zs{j}\in \cX_\zs{0}
\right)\ge (1-\delta)^{j}\,.
$$
Therefore, 
taking into account that 
for any $\check{x}\in\check{X}$ we have
$\check{\P}^{\vartheta}(\check{x},\cX_\zs{1}\setminus C_\zs{1})=0$, we obtain 
for $j\ge 2$
\begin{align*}
\check{\P}^{\vartheta}_\zs{\alpha}(\check{\tau}_\zs{\alpha}=j)&=
\check{\P}^{\vartheta}_\zs{\alpha}
\left(
\check{\Phi}_\zs{1}\notin\alpha\,,\ldots,
\check{\Phi}_\zs{j-1}\notin\alpha\,,
\check{\Phi}_\zs{j}\in\alpha
\right)
\\[2mm]
&=
\check{\P}^{\vartheta}_\zs{\alpha}
\left(
\check{\Phi}_\zs{1}\in \cX_\zs{0}\,,\ldots,
\check{\Phi}_\zs{j-1}\in \cX_\zs{0}\,,
\check{\Phi}_\zs{j}\in C_\zs{1}
\right)\\[2mm]
&\ge
\delta\,\eta_\zs{1}\,
\check{\P}^{\vartheta}_\zs{\alpha}
\left(
\check{\Phi}_\zs{1}\in \cX_\zs{0}\,,\ldots,
\check{\Phi}_\zs{j-1}\in \cX_\zs{0}\right)\,.
\end{align*}
This yields the lower bound
\eqref{sec:Pr.11}. 
Similarly we can show
$$
\sup_\zs{j\ge 1}
\frac{\sup_\zs{\vartheta\in\Theta}\,
\check{\P}^{\vartheta}_\zs{\alpha}(\check{\tau}_\zs{\alpha}=j)
}{\left(1
-\delta\eta_\zs{1}
\right)^{j-1}}
\le 
\delta
$$
and 
$$
\sup_\zs{\vartheta\in\Theta}\,
\check{\E}^{\vartheta}_\zs{\alpha} e^{-r_\zs{*}Y_\zs{1}}\le 
\frac{\delta}{e^{r_\zs{*}}-1+\delta\eta_\zs{1}}
\,.
$$
Moreover, taking into account  that 
$\check{\E}^{\vartheta}_\zs{\alpha} e^{-r_\zs{*}Y_\zs{1}}\le e^{-r_\zs{*}}$ we obtain that
$$
\sup_\zs{\vartheta\in\Theta}\,
\check{\E}^{\vartheta}_\zs{\alpha} e^{-r_\zs{*}Y_\zs{1}}\le 
\,
\check{B}^{*}_\zs{1}<1
\,,
$$
where $\check{B}^{*}_\zs{1}$ is given in \eqref{sec:De.1-1}.
Therefore, the sequence $(Y_\zs{j})_\zs{j\ge 0}$ satisfies the condition of
Theorem~\ref{Th.sec:Cr.1} with $r=r_\zs{*}$, $\upsilon_\zs{*}(r)=\ln \left(V^{*} B^{*}\right)$,
$l_\zs{*}=\wt{l}$,
$$
\epsilon_\zs{*}=\delta\,\eta_\zs{1}\,(1-\delta)^{\wt{l}-2}
\quad\mbox{and}\quad 
\upsilon_\zs{1}(r)=\ln \check{B}^{*}_\zs{1}
\,,  
$$
where $\wt{l}$ is defined in \eqref{sec:De.1-1}.
Therefore, for $n\ge 2$,
$$
\omega(n)\le \wt{A}_\zs{3}\,e^{-\wt{\kappa} n}\,,
$$
where the parameters
$\wt{A}_\zs{3}$ and
$\wt{\kappa}$ are defined
in \eqref{sec:De.5}. Taking into account that
$\wt{A}_\zs{3}\ge e^{\wt{\kappa}}$ we obtain that, for any $n\ge 0$,
$$
\omega(n)\le \wt{A}_\zs{3}\,e^{-\wt{\kappa} n}\,.
$$
Therefore, for any $0<\varkappa<\wt{\kappa}$
\begin{equation}\label{sec:Pr.12}
\wh{\omega}(\varkappa)=\sum_\zs{n\ge 0}\,
e^{\varkappa n}\omega(n)
\le \wt{A}_\zs{3}
\frac{e^{\wt{\kappa}-\varkappa}}{e^{\wt{\kappa}-\varkappa}-1}
\,.
\end{equation}
Now, taking into account that
$$
\check{\pi}^{\vartheta}(f)=
\frac{1}{\check{\E}_\zs{\alpha} Y_\zs{1}}
\,\sum^{+\infty}_\zs{j=0}\,g_\zs{f,\alpha}(j)
$$
and that $\check{\E}_\zs{\alpha} Y_\zs{1}\ge 1$,
one obtains, for any $n\ge 1$,
$$
|
\check{\E}^{\vartheta}_\zs{\check{x}}\,f(\check{\Phi}_\zs{n})\,
\Chi_\zs{\{\check{\tau}_\zs{\alpha}<n\}}
-
\check{\pi}^{\vartheta}(f)
|\,\le\,
\a_\zs{\check{x}}*\omega*g_\zs{f,\alpha}(n)
+
Q_\zs{\check{x}}*g_\zs{f,\alpha}(n)
+
G_\zs{f,\alpha}(n)\,,
$$
where
$$
Q_\zs{\check{x}}(n)=\check{\P}^{\vartheta}_\zs{\check{x}}
\left(
\tau_\zs{\alpha}>n
\right)
\quad\mbox{and}\quad
G_\zs{f,\alpha}(n)=
\sum_\zs{j=n+1}\,g_\zs{f,\alpha}(j)\,.
$$
Therefore, for any $n\ge 0$,
\begin{align}\nonumber
 \Delta_\zs{\check{x}}(n)=|
\check{\E}^{\vartheta}_\zs{\check{x}}\,f(\check{\Phi}_\zs{n})\,
-
\check{\pi}^{\vartheta}(f)
|&\le
\a_\zs{\check{x}}*\omega*g_\zs{f,\alpha}(n)\\[2mm]\label{sec:Pr.13}
&+
Q_\zs{\check{x}}*g_\zs{f,\alpha}(n)
+g_\zs{f,\check{x}}(n)
+ G_\zs{f,\alpha}(n)\,,
\end{align}
where
$g_\zs{f,\check{x}}(n)=\check{\E}^{\vartheta}_\zs{\check{x}}\,f(\check{\Phi}_\zs{n})\,
\Chi_\zs{\{\check{\tau}_\zs{\alpha}\ge n\}}$.
Now for any sequence $(b(n))_\zs{n\ge 0}$ we denote by
$\wh{b}(\varkappa)$ the Laplace transformation, i.e.
$$
\wh{b}(\varkappa)=\sum_\zs{n\ge 0}\,e^{\varkappa n}\,b(n)\,.
$$
Therefore, from \eqref{sec:Pr.13} we obtain
for any
$0<\varkappa< r_\zs{*}$,
\begin{equation}\label{sec:Pr.14}
\wh{\Delta}_\zs{\check{x}}(\varkappa)
\le \wh{a}_\zs{\check{x}}(\varkappa)
\wh{\omega}(\varkappa)\wh{g}_\zs{f,\alpha}(\varkappa)
+
\wh{Q}_\zs{\check{x}}(\varkappa)
\wh{g}_\zs{f,\alpha}(\varkappa)
+
\wh{g}_\zs{f,\check{x}}(\varkappa)
+
\wh{G}_\zs{f,\alpha}(\varkappa)
\,.
\end{equation}
It is clear 
$$
\wh{g}_\zs{f,\check{x}}(\varkappa)=
 \check{U}^{\vartheta}_\zs{\alpha}(\check{x},\varkappa,\check{f})
\quad\mbox{with}\quad
\check{f}(\check{x})=f(<x>_\zs{1})\,.
$$
Moreover, note that 
$$\wh{a}_\zs{\check{x}}(\varkappa)=\check{\E}_\zs{\check{x}}\,
e^{\varkappa\check{\tau}_\zs{\alpha}}\le \check{U}^{\vartheta}_\zs{\alpha}(\check{x},\varkappa,\check{f})
\,,
\quad
\wh{Q}_\zs{\check{x}}(\varkappa)=\frac{\check{\E}_\zs{\check{x}}\,
e^{\varkappa\check{\tau}_\zs{\alpha}}-1}{e^{\varkappa}-1}\le 
\frac{\check{U}^{\vartheta}_\zs{\alpha}(\check{x},\varkappa,\check{f})}{e^{\varkappa}-1}
$$
and
$$
\wh{G}_\zs{f,\alpha}(\varkappa)=
\frac{\wh{g}_\zs{f,\check{x}}(\varkappa)
-
\wh{g}_\zs{f,\check{x}}(0)}{e^{\varkappa}-1}
\le
\frac{\check{U}^{\vartheta}_\zs{\alpha}(\alpha,\varkappa,\check{f})}{e^{\varkappa}-1}\,.
$$
Therefore, taking into account that 
$\wt{\kappa}\le r_\zs{*}$, we obtain, that for any function $1\le f\le V$
and for any $0<\varkappa<\wt{\kappa}$,
\begin{align*}
\wh{\Delta}_\zs{\check{x}}(\varkappa)&\le
\left(
\wh{\omega}(\varkappa)+\frac{e^{\varkappa}+1}{e^{\varkappa}-1}
\right)\,\check{U}^{\vartheta}_\zs{\alpha}(\check{x},\varkappa,\check{V})
\check{U}^{\vartheta}_\zs{\alpha}(\alpha,\varkappa,\check{V})\\[2mm]
&\le
\left(
\wt{A}_\zs{3}
\frac{e^{\wt{\kappa}-\varkappa}}{e^{\wt{\kappa}-\varkappa}-1}
+\frac{e^{\varkappa}+1}{e^{\varkappa}-1}
\right)\,\check{U}^{\vartheta}_\zs{\alpha}(\check{x},\varkappa,\check{V})
\check{U}^{\vartheta}_\zs{\alpha}(\alpha,\varkappa,\check{V})
\,.
\end{align*}
Now, putting here $\varkappa=\kappa^{*}=\wt{\kappa}/2$ and taking into account
the inequality \eqref{sec:Pr.6} we get, for any $\check{x}\in\check{\cX}$,
that
$$
\wh{\Delta}_\zs{\check{x}}(\kappa^{*})
\le
\left(
\frac{(\wt{A}_\zs{3}+1) e^{\kappa^{*}}+1}{e^{\kappa^{*}}-1}
\right) V(<x>_\zs{1})\,V^{*}\,\left(B^{*}\right)^{2}
\,.
$$
Moreover, note that the chain $(\Phi_\zs{n})_\zs{n\ge 1}$ is ergodic
with the invariant measure $\pi^{\vartheta}$ defined in
\eqref{sec:Pr.7} and \eqref{subsec:AMc.14}.
By applying Proposition~\ref{Pr.subsec:AMc.5} with
$\lambda$ equals to the Dirac measure at $x$, we obtain that,
for any function $0<f\le V$,
\begin{align*}
\E^{\vartheta}_\zs{x}f(\Phi_\zs{n})-\pi^{\vartheta}(f)&=
(1-\delta)
\left(
\check{\E}^{\vartheta}_\zs{x_\zs{0}}\,\check{f}(\check{\Phi}_\zs{n})
-
\check{\pi}^{\vartheta}(\check{f})
\right)\,\Chi_\zs{\{x\in C\}}
\\[2mm]
&+
\delta
\left(
\check{\E}^{\vartheta}_\zs{x_\zs{1}}\,
\check{f}(\check{\Phi}_\zs{n})
-
\check{\pi}^{\vartheta}(\check{f})
\right)
\Chi_\zs{\{x\in C\}}
\\[2mm]
&+
\left(
\check{\E}^{\vartheta}_\zs{x_\zs{0}}\,\check{f}(\check{\Phi}_\zs{n})
-
\check{\pi}^{\vartheta}(\check{f})
\right)
\Chi_\zs{\{x\in C^{c}\}}\,.
\end{align*}
Therefore, for any $x\in \cX$, one gets
$$
|\E^{\vartheta}_\zs{x}f(\Phi_\zs{n})-\pi^{\vartheta}(f)|
\le
\,
|
\check{\E}^{\vartheta}_\zs{x_\zs{0}}\,\check{f}(\check{\Phi}_\zs{n})
-
\check{\pi}^{\vartheta}(\check{f})
|
+
|
\check{\E}^{\vartheta}_\zs{x_\zs{1}}\,
\check{f}(\check{\Phi}_\zs{n})
-
\check{\pi}^{\vartheta}(\check{f})
|
\,,
$$
i.e.
$$
\sum_\zs{n\ge 0}\,e^{\kappa^{*} n}\,
|\E^{\vartheta}_\zs{x}f(\Phi_\zs{n})-\pi^{\vartheta}(f)|
\le
\wh{\Delta}_\zs{x_\zs{0}}(\kappa^{*})
+
\wh{\Delta}_\zs{x_\zs{1}}(\kappa^{*})\,.
$$
From here it follows the inequality
\eqref{sec:Mr.3}. Hence Theorem~\ref{Th.sec:Mr.1}.
\endproof

\medskip

\section{Application to diffusion processes}\label{sec:Ug}

In order to study geometric ergodicity for the process (\ref{sec:In.3}) we  start with
the chain $(\Phi^y_n)_\zs{n\ge 0}$, where $\Phi^y_n\,=\,y_n$.

\begin{proposition}\label{Pr.sec:Ug.1}
For any $\vartheta\in\Theta$, the sequence $(\Phi_n^y)_\zs{n\ge 0}$ is a homogeneous Markov chain aperiodic and $\psi$-irreducible, where $\psi$ is the Lebesgue measure
on $\cB(\bbr)$.
\end{proposition}
\proof
Taking into account (see, for example, \cite{GiSk}) that the solution of the equation
\eqref{sec:In.3} is a homogeneous Markov process, we obtain immediately that $(\Phi_n^y)_\zs{n\ge 0}$ is a homogeneous Markov chain. In this case (see \cite{GiSk}), for any $\vartheta$, the process $(y_\zs{t})_\zs{t\ge 0}$ admits
the transition density $\upsilon_\zs{\vartheta}(t,x,y)$ as follows :
\begin{equation}\label{sec:Ug.1}
\upsilon_\zs{\vartheta}(t,x,y)=\frac{\Upsilon(t,x,y)}{\sqrt{2\pi t}\sigma(y)}
\exp\left\{\int_{\varsigma(x)}^{\varsigma(y)}\,H_\zs{\vartheta}(u)\d u-
\frac{(\varsigma(y)-\varsigma(x))^2}{2t}\right\}\,,
\end{equation}
Here,
$$
H_\zs{\vartheta}(z)=\frac{S(\check{\varsigma}(z))}{\sigma(\check{\varsigma}(z))}-
\frac{\dot{\sigma}(\check{\varsigma}(z))}{2 \sigma^{2}(\check{\varsigma}(z)}
\,,\quad
\varsigma(x)=\int_\zs{0}^{x}\,\sigma^{-1}(u)\,\d u
$$
and $\check{\varsigma}(\cdot)$ is the inverse function of
$\varsigma$, i.e. it
is the unique solution of the equation
$z=\varsigma(\check{\varsigma})$. Moreover,
$$\Upsilon(t,x,y)=\E\,\exp\left\{-\frac{1}{2}\int_\zs{0}^t\,
\wt{H}_\zs{\vartheta}(w^{*}_\zs{u,t}(x,y))\d u\right\}\,,
$$
$\wt{H}_\zs{\vartheta}(x)=\dot{H}_\zs{\vartheta}(x)+H^2_\zs{\vartheta}(x)$
and
$$
w^{*}_\zs{u,t}(x,y)=\varsigma(x)+
\frac{u}{t}(\varsigma(y)-\varsigma(x))+w_\zs{u}-\frac{u}{t}\,w_\zs{t}\,.
$$
It  means that, for any $n\,\ge\,1$, for any $A\,\in\,\cB(\bbr)$ and for any $x\in\bbr$,
\begin{equation}\label{sec:Ug.2}
\P^{\vartheta}(\Phi_n^y\,\in\,A|\Phi_\zs{0}^{y}=x)=\int_\zs{A}
\upsilon_\zs{\vartheta}(n,x,z)\d z\,.
\end{equation}
Thus, the chain $(\Phi_n^y)_\zs{n\ge 0}$ is $\psi$-irreducible, where
$\psi$ is the  Lebesgue measure on $\cB(\bbr)$. Moreover, in this case the
chain is aperiodic.
\endproof

\noindent Now, we check the minorization condition
$\H_\zs{1})$ for the chain $(\Phi^y_\zs{n})_\zs{n\ge 0}$.

\medskip

\begin{proposition}\label{Pr.sec:Ug.2}
For any $K\ge 8\sigma_\zs{min}$,
the chain $(\Phi_\zs{n}^{y})_\zs{n\ge 0}$ satisfies the minorization condition
$\H_\zs{1})$ with $C=[-K, K]$, $\delta=\delta_\zs{K}$ and $\eta=\eta_\zs{K}$
defined in \eqref{sec:De.8}, and
the probability measure $\nu_\zs{K}$ defined in \eqref{sec:De.6-0}.
\end{proposition}
\proof
We start with studying the properties of the function $H_\zs{\vartheta}$ defined in
\eqref{sec:Ug.1}. From definition of the class $\Theta$ we find immediately that,
for any $z\in\bbr$,
$$
|\check{\varsigma}(z)|\le |z|\,\sigma_\zs{max}\,.
$$
Moreover, note that, for any $\vartheta$ from $\Theta$ and for any $y$ from $\bbr$,
\begin{equation}\label{sec:Ug.3-1}
|\dot{S}(y)|\le M+L
\quad\mbox{and}\quad
|S(y)|\le M+L|y|\,.
\end{equation}
Therefore,
$$
\sup_\zs{|z|\le z_\zs{*}}\,\sup_\zs{\vartheta\in\Theta}
\,
|H_\zs{\vartheta}(z)|
\le
H^{*}_\zs{0}
+
L\sigma_\zs{*}
\,
 z_\zs{*}\,,
$$
where  $H^{*}_\zs{0}$ is given in \eqref{sec:De.6-0}.
Note now, that the derivative of $H_\zs{\vartheta}$ can be represented as
$\dot{H}_\zs{\vartheta}(z)=F_\zs{\vartheta}(\check{\varsigma}(z))$ with
$$
F_\zs{\vartheta}(y)=\dot{S}(y)-\frac{S(y)\dot{\sigma}(y)}{\sigma(y)}
-\frac{\ddot{\sigma}(y)}{2\sigma(y)}+
\frac{(\dot{\sigma}(y))^{2}}{\sigma^{2}(y)}\,.
$$
By making use of the upper bounds \eqref{sec:Ug.3-1}, we obtain
$$
\sup_\zs{\vartheta\in\Theta}
|F_\zs{\vartheta}(y)|\le H^{*}_\zs{1}+ L\sigma_\zs{*} |y|
$$
and, therefore, for any $z_\zs{*}>0$,
$$
\sup_\zs{|z|\le z_\zs{*}}\,\sup_\zs{\vartheta\in\Theta}
|\wt{H}_\zs{\vartheta}(z)|
\le \wt{H}_\zs{*}(z_\zs{*})
\,,
$$
where 
$\wt{H}_\zs{*}(z)=H^{*}_\zs{1}+2(H^{*}_\zs{0})^{2}+
L\sigma_\zs{*}\sigma_\zs{max}\,z
+2(L\sigma_\zs{*})^2\,
z^{2}$
and the coefficient 
$H^{*}_\zs{1}$ is given in \eqref{sec:De.6-0}.
Moreover, we note that on the set
$$
\Gamma_\zs{K}=
\{\sup_\zs{0\le u\le 1}|w_\zs{u}|\le \wt{K}/2\}\quad\mbox{with}\quad \wt{K}=K/\sigma_\zs{min}\,,
$$
the process $(w^{*}_\zs{v,t}(x,y))_\zs{0\le v\le t\le 1}$  is bounded :
$$
\sup_\zs{x,y\in C}\,
\sup_\zs{0\le u\le t\le 1}
|w^{*}_\zs{u,t}(x,y)|\le 2\wt{K}\,.
$$
Therefore, for any $x,y$ from $C$,
$$
\Upsilon(1,x,y)\ge \P(\Gamma_\zs{K})\,
e^{-\frac{1}{2}\wt{H}_\zs{*}(2\wt{K})}\,.
$$
By making use of the Doob inequality we obtain
$$
\P(\Gamma_\zs{K})\ge 1-\frac{4\E\sup_\zs{0\le t\le 1}w^2_\zs{t}}{K^{2}_\zs{1}}
\ge 1-\frac{16 \sigma^{2}_\zs{min}}{K^{2}}\,,
$$
i.e. for $K\ge 8 \sigma_\zs{min}$,
$$
\P(\Gamma_\zs{K})\ge 3/4
\quad\mbox{and}\quad
\Upsilon(1,x,y)\ge \frac{3}{4}\,
e^{-\frac{1}{2}\wt{H}_\zs{*}(2\wt{K})}\,.
$$
Moreover, for any $x,y\in C$,
$$
\sup_\zs{\vartheta\in\Theta}\,
\left|
\int_{\varsigma(x)}^{\varsigma(y)}\,H_\zs{\vartheta}(u)\d u
\right|
\le 
\,
\left( H^{*}_\zs{0}+L\sigma_\zs{*} \wt{K} 
\right)\, \wt{K}\,.
$$
This implies that
$$
\inf_\zs{x,z\in C}\,
\upsilon_\zs{\vartheta}(1,x,z)\ge
3(4\sqrt{2\pi} \sigma_\zs{max})^{-1}\,e^{-\Omega_\zs{*}(\wt{K})}\,,
$$
where $\Omega_\zs{*}(z)$ is introduced in \eqref{sec:De.8}.
Therefore, taking into account that, for any
$A$ from $\cB(\cX)$,
$$
\P^{\vartheta}(\Phi_1^y\,\in\,A|\Phi^y_\zs{0}=x)=\int_A\,\upsilon_\zs{\vartheta}(1,x,z)\d z\,,
$$
yields the inequality \eqref{sec:Mr.2} with
$\delta_\zs{K}$ and $\nu_\zs{K}(\cdot)$ defined
in \eqref{sec:De.8} and \eqref{sec:De.6-0} respectively.
Moreover, the inequality \eqref{sec:Mr.2-1} follows 
directly from Proposition~\ref{Pr.subsec:AMc.8}.
Hence Proposition~\ref{Pr.sec:Ug.2}.
\endproof

Now, for any $\C^{2}(\bbr)\to \bbr$ function $V$, we introduce the generator
\begin{equation}\label{sec:Ug.5}
\A_\zs{\vartheta}(V)(x)=
\dot{V}(x)S(x)+
\frac{1}{2}\,
\sigma^{2}(x)\ddot{V}(x)\,.
\end{equation}

\begin{definition}\label{De.sec:Ug.1}
Any $\bbr\to [1,\infty)$ twice
continuously differentiable function $V$ is called uniform
over $\vartheta\in\Theta$ Lyapunov function for the equation \eqref{sec:In.3} if
the following conditions fulfill:

\begin{itemize}
 \item 
for some constants $\gamma>0$, $b^{*}>0$ and for any $x\in\bbr$,
\begin{equation}\label{sec:Ug.6}
\sup_\zs{\vartheta\in\Theta}
\A_\zs{\vartheta}(V)(x)
\le -\gamma V(x)+b^{*}\,;
\end{equation}

\item
$\lim_{x\to\infty}V(x)=\infty$ and there exists $m>0$ such that
\begin{equation}\label{sec:Ug.7}
\sup_\zs{x\in\bbr}\frac{V(x)+|\dot{V}(x)|}{1+|x|^m} < \infty\,.
\end{equation}
\end{itemize}
\end{definition}

\medskip

\begin{proposition}\label{Pr.sec:Ug.2-1}
For any $0<\epsilon\le 1/2$, the function
$V(x)=(1+x^2)^{\epsilon}$
satisfies the inequality \eqref{sec:Ug.6}
with $\gamma= 2\epsilon \beta$ and $b^{*}=\epsilon b^{*}_\zs{0}$,
where $b^{*}_\zs{0}$ is given in \eqref{sec:De.8-1}.
\end{proposition}
\proof
The definition of the space $\Theta$ implies directly that,
for $|x|\ge \x_\zs{*}$,
$$
xS(x)\le |x|(M+\beta\x_\zs{*})-\beta x^{2}\,.
$$
Therefore, we get, for any $\vartheta\in\Theta$,
\begin{align*}
\A_\zs{\vartheta}(V)(x)
&\le
\frac{2\epsilon V(x)xS(x)}{1+x^2}
+\epsilon \sigma^{2}_\zs{max}
\\[2mm]
&\le \frac{2\epsilon\,V(x)xS(x)}{1+x^2}\Chi_\zs{\{|x|\ge \x_\zs{*}\}}
+\epsilon
\left(
(1+\x^{2}_\zs{*})^{\epsilon}M+\sigma^{2}_\zs{max}
 \right)
\\[2mm]
&\le -2\epsilon V(x)\beta +b^{*}\,.
\end{align*}
Hence Proposition~\ref{Pr.sec:Ug.2-1}.
\endproof

\medskip
\medskip
\medskip
\begin{proposition}\label{Pr.sec:Ug.3}
 Let $V(x)$ be a uniform over $\vartheta\in\Theta$ Lyapunov function for equation
\eqref{sec:In.3} from definition~\ref{De.sec:Ug.1} with constants $\gamma$ and $b^{*}$. Then, for 
any $K>0$ and $0<\check{\varepsilon}<1$ for which
$$
\inf_\zs{|x|\ge K}\,V(x)
\ge
\,\frac{b^{*}}{\check{\varepsilon} \gamma}\,,
$$
the chain
$(\Phi_n^y)_\zs{n\ge 0}$ satisfies the following inequality
\begin{equation}\label{sec:Ug.9}
\sup_\zs{\vartheta\in\Theta}\,\E^{\vartheta}_\zs{x}\,V(\Phi_1^y)\le\,
(1-\rho) V(x)+D_\zs{K}
\,
\Chi_\zs{\{|x|\le K\}}\,,
\end{equation}
where $\rho=(1-\check{\varepsilon})(1-e^{-\gamma})$,
$D_\zs{K}=V^{*}_\zs{K}e^{-\gamma}+b^{*}(1-e^{-\gamma})/\gamma$
 and
$V^{*}_\zs{K}=\sup_\zs{|x|\le K}\,V(x)$.
\end{proposition}
\proof
By the Ito formula, one gets
\begin{align*}
V(y_\zs{t})&=V(y_\zs{0})+\int^{t}_\zs{0}\,
\A_\zs{\vartheta}(V)(y_\zs{s})\d s+\,\int^{t}_\zs{0}\,\dot{V}(y_\zs{s})\sigma(y_\zs{s})\d w_\zs{s}\,.
\end{align*}
In Proposition~\ref{Pr.subsec:AMc.7}, we have proved that the moments of the solution
of equation (\ref{sec:In.3}) are bounded. This implies that the stochastic integral
is a martingale in the above Ito formula. Therefore, by setting
$Z(t)=\E^{\vartheta}_\zs{x}\,V(y_\zs{t})$,
 one has
$$
\dot{Z}(t)=
\E^{\vartheta}_\zs{x}
\,\A_\zs{\vartheta}(V)(y_\zs{t})=-\gamma Z(t)+\psi_\zs{t}
\,,
$$
where $\psi_\zs{t}=\E^{\vartheta}_\zs{x}\left(\A_\zs{\vartheta}(V)(y_\zs{t})+\gamma y_\zs{t}\right)$.
The inequality \eqref{sec:Ug.6} gives $\psi_\zs{t}\le b^{*}$. Resolving this differential equation,
we obtain, that for $0\le t\le 1$
\begin{align*}
Z(t)&= Z(0)e^{-\gamma t}+
\int^{t}_\zs{0}\,e^{-\gamma(t-s)}\psi_\zs{s}\d s\\[2mm]
&\le  Z(0)e^{-\gamma t}+
b^{*}\,
\frac{1-e^{-\gamma}}{\gamma}
=V(x)e^{-\gamma t}+b^{*}\frac{1-e^{-\gamma}}{\gamma}\,.
\end{align*}
Therefore,
$$
\sup_\zs{\vartheta\in\Theta}\,
\E^{\vartheta}_\zs{x}\,V(\Phi_1^y)\,
\le V(x)e^{-\gamma }+
b^{*}
\frac{1-e^{-\gamma}}{\gamma}
\,.
$$
From here we obtain the inequality \eqref{sec:Ug.9}.
Hence Proposition~\ref{Pr.sec:Ug.3}. \endproof

\subsection{Proof of Theorem~\ref{Th.sec:Mr.2}}

First note, that thanks to Propositions~\ref{Pr.sec:Ug.1}
 the diffusion process \eqref{sec:In.3} with the parameter
$\vartheta=(S,\sigma)$ from $\Theta$ introduced in  \eqref{sec:Mr.3-1},
 satisfies the condition
$\H_\zs{1})$ with the set $C=[-K,K]$ for $K \ge K_\zs{0}$ given \eqref{sec:De.6-1}
 and the measure $\nu$ defined
in \eqref{sec:De.6-0}. Moreover, Propositions~\ref{Pr.sec:Ug.2-1} -- \ref{Pr.sec:Ug.3}
imply that for any $0<\epsilon\le 1/2$ this process satisfies the condition
$\H_\zs{2})$ with the parameters \eqref{sec:De.9}.
Moreover, for any $t\ge 1$ and any $\bbr\to ]0,1]$ function $g$, we set
$$
\wt{g}(x)=\E^{\vartheta}_\zs{x}\,g(y_\zs{t})\,=
\E^{\vartheta}_\zs{x}\,g(y_\zs{\{t\}})\,.
$$
Moreover, taking into account that $\pi(g)=\pi(\tilde{g})$, one has
$$
\E^{\vartheta}_\zs{x}\,g(y_\zs{t})-\pi(g)=
\E^{\vartheta}_\zs{x}\,\wt{g}(\Phi_\zs{[t]}^{y})-\pi(\tilde{g})\,.
$$
Therefore by applying Theorem~\ref{Th.sec:Mr.1}
 to the chain $(\Phi_n^y)_\zs{n\ge 0}$, we come to the upper bound
\eqref{sec:Mr.4}. Hence Theorem~\ref{Th.sec:Mr.2}.
 \endproof

\medskip
\medskip

\setcounter{section}{0}
\renewcommand{\thesection}{\Alph{section}}

\section{Appendix}\label{sec:A}

\subsection{Homogeneous Markov chains with atoms}\label{subsec:AMc}
We follow the Meyn-Tweedie approach (see \cite{MeTw}). We remind some definitions from \cite{MeTw} for a
 homogeneous Markov chains $(\Phi_n)_\zs{n\ge 0}$ defined on a measurable state space
$(\cX, \cB(\cX))$. Denote by $P(x,\cdot)\,, x\in\cX\,,$ the transition probability of this chain, i.e.
for any $A\in\cB(\cX), x\in\cX$,
$$
P(x,A)\,=\,\P_\zs{x}(\Phi_\zs{1}\in A)
=
\P(\Phi_1\in A|\Phi_0\,=\,x)\,.
$$
The $n-$step transition probability is
$$
P^{n}(x,A)\,=\,\P_\zs{x}(\Phi_n\in A)\,.
$$
We remind that a measure $\pi$ on $\cB(\cX))$ is called {\it invariant} for this chain if,
for any $A\in\cB(\cX)$,
$$
\pi(A)\,=\,\int_{\cX}\,P(x,A)\pi(\d x)\,.
$$

If there exists an invariant positive measure $\pi$ with $\pi(\cX)\,=\,1$ then the chain is
called {\it positive}.
\begin{definition}\label{De.subsec:AMc.1}
The chain $(\Phi_n)_\zs{n\ge 0}$ is $\varphi$-irreducible if there exists a
nontrivial  measure $\varphi$
on $\cB(\cX)$ such that, whenever $\varphi(A)\,>\,0$, one has
$$
L(x,A)=\P_\zs{x}(\cup_\zs{n=1}^{\infty}\{\Phi_n\,\in\,A\})\,>\,0 \quad\hbox{\rm for any}\quad x\,\in\,\cX\,.
$$
\end{definition}
One can show that, for any $\varphi$-irreducible chain, there exists a "maximal" irreducible measure which is noted as $\psi$ and the chain is called $\psi$-irreducible. A irreducible measure
$\psi$ is maximal if and only if  $\psi(A)=0$ implies
$$
\psi(x\in\bbr\,:\,L(x,A)>0)=0\,.
$$
In the sequel, we denote
$$
\cB_\zs{+}(\cX)\,=\,\{A\,\in\,\cB(\cX) : \psi(A)\,>\,0\}\,.
$$
\begin{definition}\label{De.subsec:AMc.2}
The chain $(\Phi_n)_\zs{n\ge 0}$ is Harris recurrent if it is $\psi$-irreducible and, for
any $A\,\in\,\cB_+(\cX)$,  one has
$$
\P_\zs{x}\left(\sum_\zs{n=1}^{\infty}\,
\Chi_\zs{\{\Phi_n\in A\}}\right)=1\,, \quad\hbox{\rm for any}\quad x\,\in\,A\,.
$$
\end{definition}

\begin{definition}\label{De.subsec:AMc.3}
The Markov $\psi$-irreducible chain $(\Phi_n)_\zs{n\ge 0}$ is called periodic of period $d$ if
there exist disjoint sets $\Gamma_1,\ldots,\Gamma_d$ in $\cB(\cX)$ with
$$
\psi\left(\cap_\zs{j=1}^{d}\Gamma_j^c\right)\,=\,0
$$
such that, for $1\,\le\,i\,\le\,d-1$ and for any $x\,\in\,\Gamma_i$, one has
$$
\P_\zs{x}(\Phi_1\,\in\,\Gamma_\zs{i+1})\,=\,1
$$
and for $x\,\in\,\Gamma_d$ one has $\P_\zs{x}(\Phi_1\,\in\,\Gamma_\zs{1})=1$.
The chain is aperiodic if $d\,=\,1$.
\end{definition}

\begin{definition}\label{De.subsec:AMc.4}
We will say that the chain $(\Phi_n)_\zs{n\ge 0}$ satisfies the minorization condition if, for some $\delta\,>\,0$, some set $C\,\in\,\cB(\cX)$ and some probability measure $\nu$ with
$\nu(C)\,=\,1$, one has
\begin{equation}\label{subsec:AMc.1}
\inf_\zs{A\in \cB(\cX)}
\left(
\inf_\zs{x\in C}
P(x,A)-\delta\,\nu(A)\right)
\ge 0\,.
\end{equation}
\end{definition}

\begin{definition}\label{De.subsec:AMc.6}
A set $\alpha\in\cB_\zs{+}(\cX)$ is called accessible atom if,
for any $x$ and $y$ from $\alpha$,
$$
\P(x,\Gamma)=\P(y,\Gamma)\,,\quad \forall\  \Gamma\in\cB(\cX)\,.
$$
 \end{definition}
In order to study the ergodicity property, we associate to any set $C\in\cB(\cX)$
the stopping time
$$
\tau_\zs{C}=\inf\{k\ge 1\,:\,\Phi_\zs{k}\in C\}\,.
$$
\begin{proposition}\label{Pr.subsec:AMc.0}
Suppose that the Markov chain $\Phi$ is $\psi$-irreducible and contains an accessible atom
$\alpha$ such that
$$
\E_\zs{\alpha}\tau_\zs{\alpha}
<\infty\,.
$$
Then the chain is ergodic with the
invariant probability measure $\pi$ defined
as
$$
\pi(\Gamma)
=
\frac{1}{\E_\zs{\alpha}\,\tau_\zs{\alpha}}\,
\E_\zs{\alpha}\,
\sum^{\tau_\zs{\alpha}}_\zs{j=1}\,\Chi_\zs{\{\Phi_\zs{j}\in\Gamma\}}\,.
$$
\end{proposition}
\proof
Indeed, by the definition of $\pi$, for any set $\Gamma\in\cB(\cX)$, one has
\begin{align*}
\int_\zs{\cX}\,\pi(\d z)\,\P(z,\Gamma)&=
\frac{1}{\E_\zs{\alpha}\,\tau_\zs{\alpha}}\,
\E_\zs{\alpha}\,
\sum^{\infty}_\zs{j=1}\,
\Chi_\zs{\{j\le\tau_\zs{\alpha}\}}
\,\E_\zs{\alpha}\left(
\Chi_\zs{\{\Phi_\zs{j+1}\in\Gamma\}}
|\Phi_\zs{1},\ldots,\Phi_\zs{j}
\right)
\\[2mm]
&=\frac{1}{\E_\zs{\alpha}\,\tau_\zs{\alpha}}\,
\left(
\E_\zs{\alpha}\,
\sum^{\tau_\zs{\alpha}}_\zs{j=2}\,\Chi_\zs{\{\Phi_\zs{j}\in\Gamma\}}
+
\P_\zs{\alpha}\,\left(
\Phi_\zs{\tau_\zs{\alpha}+1}\in\Gamma \right)
\right)
\,.
\end{align*}
Moreover, it is easy to see that
$$
\P_\zs{\alpha}
\left(
\Phi_\zs{\tau_\zs{\alpha}+1}\in\Gamma
\right)
=
\P_\zs{\alpha}
\left(
\Phi_\zs{1}\in\Gamma
\right)\,.
$$
This implies the relationship
$$
\int_\zs{\cX}\,\pi(\d z)\,\P(z,\Gamma)=\pi(\Gamma)\,,
$$
i.e. the measure $\pi$ is invariant. Obviously, that $\pi(\cX)=1$, i.e.
$\pi$ is a probability measure.
\endproof

\medskip

\medskip

\subsection{Lyapunov functions method for Markov chains}

We start with the definition of a "Lyapunov function". For this we impose the following
drift condition to the chain $(\Phi_n)_\zs{n\ge 0}$, i.e.

$\D)$ {\em
There exist
a $\cX\to [1,\infty)$ function $V$, constants $0<\rho<1$,  $D\ge 1$ and
 a set $C$ from $\cB(\cX)$ such that
for all $x\in\cX$}
\begin{equation}\label{subsec:AMc.2}
\E_\zs{x}\left(V(\Phi_\zs{1})\right)\,
\le\,(1-\rho)V(x)\,+\,D\Chi_\zs{C}(x)\,.
\end{equation}
In this case we call $V$ the {\em Lyapunov function.}


Now, for any $\cX\to [1,+\infty)$ function $f$ and any set $A\in\cB(\cX)$,
we set
\begin{equation}\label{subsec:AMc.2-1}
U_\zs{A}(x,r,f)=
\E_\zs{x}\,\sum^{\tau_\zs{A}}_\zs{j=1}\,e^{r j}\,f(\Phi_\zs{j})\,.
\end{equation}

\begin{proposition}\label{Pr.subsec:AMc.1}
Assume, that  
the condition $\D)$ holds. 
 Then,
for any\\ $0<a<-\ln(1-\rho)$, one has
\begin{equation}\label{subsec:AMc.2-2}
\sup_\zs{x\in\cX}\,
\frac{U_\zs{C}(x,a,V)}{V(x)}
\le
\,U^{*}(a)\,,
\end{equation}
where
$$
U^{*}(a)=\frac{(1-\rho)e^{a}+D\,e^{a}}{1-(1-\rho)e^{a}}
\,.
$$
\end{proposition}
\proof
The condition \eqref{subsec:AMc.2} implies immediately
$$
U_\zs{C}(x,a,V)\le
\,(1-\rho)e^{a}V(x)+
\,(1-\rho)e^{a}U_\zs{C}(x,a,V)+
D\,e^{a}\,.
$$
Taking into account that $V(x)\ge 1$, we obtain
the inequality \eqref{subsec:AMc.2-2}.
\endproof

\begin{proposition}\label{Pr.subsec:AMc.3}
Assume that for some $a=a_\zs{*}>0$ the Markov chain $(\Phi_\zs{n})_\zs{n\ge 1}$
satisfies the property \eqref{subsec:AMc.2-2}
with the $C$ - bounded function
$V$,  i.e.
\begin{equation}\label{subsec:AMc.8}
V^{*}=\sup_\zs{x\in C}\,V(x)<\infty\,.
\end{equation}
Let $B$  be a set from $\cB(\cX)$ such that, for some $\k_\zs{*}>0$,
\begin{equation}\label{subsec:AMc.8-1}
\inf_\zs{x\in C\setminus B}\,\P(x,B)\ge \k_\zs{*}\,.
\end{equation}
Then, for any $0<\gamma<1$ and
\begin{equation}\label{subsec:AMc.8-2}
0<a\le 
\,
\frac{(1-\gamma)|\ln(1-\k_\zs{*})|}{\ln V^{*} U^{*}(a_\zs{*})+|\ln(1-\k_\zs{*})|}
\,a_\zs{*}
\,,
\end{equation}
one has
\begin{equation}\label{subsec:AMc.9}
\sup_\zs{x\in\cX}\,\frac{1}{V(x)}\,
U_\zs{B}(x,a,V)
\le U^{*}_\zs{1}(a)\,,
\end{equation}
where
$$
U^{*}_\zs{1}(a)=U^{*}(a)
\left(
1+
\frac{U^{*}(a)V^{*}}{1-(1-\k_\zs{*})^{\gamma}}
\right)
\,.
$$
\end{proposition}
\proof
First, we introduce the sequence
of  stopping times $(\tau_\zs{C}(n))_\zs{n\ge 0}$ as follows : 
$\tau_\zs{C}(0)=0$ and, for $n\ge 1$,
$$
\tau_\zs{C}(n)=\inf\{k\ge \tau_\zs{C}(n-1)+1\,:\,\Phi_\zs{k}\in C\}\,.
$$
Obviously, that $\tau_\zs{C}(1)=\tau_\zs{C}$. Moreover, the condition \eqref{subsec:AMc.2-2} 
implies that $\E_\zs{x}\tau_\zs{C}(n)<\infty$ for any $n\ge 1$ and $x\in\cX$.
Using this sequence
we obtain that for $0<a\le a_\zs{*}$
\begin{align*}
U_\zs{B}(x,a,V)&=
\sum^{\infty}_\zs{n=0}
\E_\zs{x}\,\sum^{\tau_\zs{C}(n+1)}_\zs{j=\tau_\zs{C}(n)+1}\,
e^{a j}\,V(\Phi_\zs{j})\,
\Chi_\zs{\{\tau_\zs{B}\ge j\}}\\[2mm]
&
\le
U_\zs{C}(x,a,V)
+
\sum^{\infty}_\zs{n=1}
\E_\zs{x}\,\Chi_\zs{\{\tau_\zs{B}> \tau_\zs{C}(n)\}}\,
e^{a \tau_\zs{C}(n)}\,
U_\zs{C}(z_\zs{n},a,V)\,,
\end{align*}
where 
$z_\zs{n}=\Phi_\zs{\tau_\zs{C}(n)}$.
Moreover, taking into account the inequalities \eqref{subsec:AMc.2-2} 
and \eqref{subsec:AMc.8}, we get that,
for $a\le a_\zs{*}$,
$$
U_\zs{B}(x,a,V)\le U^{*}(a)V(x)
+
U^{*}(a)V^{*}\,
\sum^{\infty}_\zs{n=1}
\Upsilon_\zs{n}(x,a)\,,
$$
where
$$
\Upsilon_\zs{n}(x,a)=\E_\zs{x}\,
\Chi_\zs{\{\tau_\zs{B}>\tau_\zs{C}(n)\}}\,e^{a \tau_\zs{C}(n)}\,.
$$
Note now that, for $n=1$,
$$
\Upsilon_\zs{1}(x,a)\le V(x) U^{*}(a)\,.
$$
Let now $n\ge 2$. We have
\begin{align*}
\Upsilon_\zs{n}(x,a)&=\E_\zs{x}\,
\Chi_\zs{\{\tau_\zs{B}>\tau_\zs{C}(n)\}}\,e^{a \tau_\zs{C}(n)}\\[2mm]
&=\E_\zs{x}\,
\Chi_\zs{\{\tau_\zs{B}>\tau_\zs{C}(n-1)\}}\,e^{a \tau_\zs{C}(n-1)}
\E_\zs{z_\zs{n-1}}\,e^{a\tau_\zs{C}}\,\Chi_\zs{\{\tau_\zs{B}>\tau_\zs{C}\}}\,.
\end{align*}
Therefore, for $n\ge 2$,
$$
\Upsilon_\zs{n}(x,a)\le  \Upsilon_\zs{n-1}(x,a)\,
\sup_\zs{z\in C \setminus B }\,
\E_\zs{z}\,e^{a\tau_\zs{C}}\,\Chi_\zs{\{\tau_\zs{B}>1\}}\,.
$$
By the H\"older inequality and the condition 
\eqref{subsec:AMc.8}, we get
$$
\sup_\zs{z\in C \setminus B }\,
\E_\zs{z}\,e^{a\tau_\zs{C}}\,\Chi_\zs{\{\tau_\zs{B}>1\}}\,\le\,
\left(V^{*} U^{*}(a_\zs{*})\right)^{a/a_\zs{*}}\,
(1-\k_\zs{*})^{1-a/a_\zs{*}}:=\g(a)\,.
$$
Therefore, for any $a$ satisfying the condition \eqref{subsec:AMc.8-2}
we obtain
$$
\g(a)\le (1-\k_\zs{*})^{\gamma}
$$
and, for any $n\ge 2$,
$$
\Upsilon_\zs{n}(x,a)\le  \Upsilon_\zs{1}(x,a)\,\g(a)^{n-1}
\le 
V(x) U^{*}(a)\,(1-\k_\zs{*})^{\gamma(n-1)}
\,.
$$
This implies directly the bound
\eqref{subsec:AMc.9}. Hence Proposition~\ref{Pr.subsec:AMc.3}.
\endproof

\medskip

\subsection{Properties of splitting chains}

Now, we study some property of the splitting chain
$(\check{\Phi}_\zs{n})_\zs{n\ge 1}$ constructed in Section 4
, which we represent as
\begin{equation}\label{subsec:AMc.11}
\check{\Phi}_\zs{n}=(\check{\phi}_\zs{n},\check{\iota}_\zs{n})\,,
\end{equation}
where $\check{\phi}_\zs{n}\in\cX$ and $\check{\iota}_\zs{n}\in\{0,1\}$.
\begin{proposition}\label{Pr.subsec:AMc.4}
For any measure $\lambda$ on $\cB(\cX)$ and any set
$\check{\Gamma}\in\cB(\check{\cX})$,
\begin{equation}\label{subsec:AMc.12}
\int_\zs{\check{\cX}}\,
\check{\P}^{\vartheta}(\check{x},\check{\Gamma})
\,\lambda^{*}(\d \check{x})
=\lambda^{*}_\zs{1}(\check{\Gamma})
\,,
\end{equation}
where
$$
\lambda_\zs{1}(\cdot)=\int_\zs{\cX}\,
\P^{\vartheta}(x,\cdot)\,\lambda(\d x)\,.
$$
\end{proposition}
\proof Indeed, by the definition of
the $*$ operation and
of the transition probability
$\check{\P}^{\vartheta}(\cdot,\cdot)$
we obtain
\begin{align*}
\int_\zs{\check{\cX}}\,
\check{\P}^{\vartheta}(\check{x},\check{\Gamma})
\,\lambda^{*}(\d \check{x})
=
\int_\zs{\cX}\,
\P^{\vartheta}(x,\check{\Gamma})^{*}
\,\lambda(\d x)
=
\lambda^{*}_\zs{1}(\check{\Gamma})\,.
\end{align*}
\endproof

\begin{proposition}\label{Pr.subsec:AMc.5}
For any $n\ge 1$, any measurable positive
$\cX^{n}\to\bbr$ function
$G$ and for any measure $\lambda$ on $\cB(\cX)$, one has
\begin{equation}\label{subsec:AMc.13}
\int_\zs{\cX}\,
\E^{\vartheta}_\zs{x}\,G_\zs{n}(\Phi_\zs{1},\ldots,\Phi_\zs{n})
\,\lambda(\d x)
=
\int_\zs{\check{\cX}}\,
\check{\E}^{\vartheta}_\zs{\check{x}}\,
G_\zs{n}(\check{\phi}_\zs{1},\ldots,\check{\phi}_\zs{n})
\,\lambda^{*}(\d \check{x})\,.
\end{equation}
\end{proposition}
\proof It is clear, that it suffices to check this equality
for  positive functions of the form
$$
G_\zs{n}(x_\zs{1},\ldots,x_\zs{n})
=
\prod^{n}_\zs{j=1}\,g_\zs{j}(x_\zs{j})\,.
$$
First, we check this equality for $n=1$. Note that, for
any $\cX\to \bbr$ function $g$
and for any $x\in\cX$, one has
$$
\int_\zs{\check{\cX}}\,g(<\check{y}>_\zs{1})\,
\P^{*}(x,\d \check{y})
=
\int_\zs{\cX}\,g(y)\,
\P(x,\d y)\,,
$$
where $<\check{y}>_\zs{1}$ denotes the first component of the
$\check{y}\in\check{\cX}=\cX\times \{0,1\}$. Making use of this equality
implies easy  \eqref{subsec:AMc.13} for $n=1$. Assume now
that the equality  \eqref{subsec:AMc.13} is true until $n-1$.
We check it for $n$. Indeed, we have
$$
\check{\E}^{\vartheta}_\zs{\check{x}}
G_\zs{n}(\check{\phi}_\zs{1},\ldots,\check{\phi}_\zs{n})
=
\check{\E}^{\vartheta}_\zs{\check{x}}\,
\prod^{n}_\zs{j=1}\,g_\zs{j}(x_\zs{j})
=\check{\E}^{\vartheta}_\zs{\check{x}}\,g_\zs{1}(\check{\phi}_\zs{1})\,
T(\check{\Phi}_\zs{1})\,,
$$
where
$$
T(\check{y})=
\check{\E}^{\vartheta}_\zs{\check{y}}\,\prod^{n-1}_\zs{j=1}\,
g_\zs{j+1}(\check{\phi}_\zs{j})\,.
$$
Now, we set
$$
\mu(\Gamma)=\int_\zs{\Gamma}\,g_\zs{1}(y)\,\lambda_\zs{1}(\d y)\,,
$$
where the measure $\lambda_\zs{1}(\cdot)$
 is defined in \eqref{subsec:AMc.12}.
Therefore, taking into account Proposition~\ref{Pr.subsec:AMc.4}
, we can represent the integral on the right hand side
of the equality \eqref{subsec:AMc.13} as
$$
\int_\zs{\check{\cX}}\,
\check{\E}^{\vartheta}_\zs{\check{x}}\,
G_\zs{n}(\check{\phi}_\zs{1},\ldots,\check{\phi}_\zs{n})
\,\lambda^{*}(\d \check{x})
=
\int_\zs{\check{\cX}}\,T(\check{y})\mu^{*}(\d \check{y})
\,.
$$
By the induction assumption, one has
$$
\int_\zs{\check{\cX}}\,T(\check{y})\mu^{*}(\d \check{y})
=
\int_\zs{\cX}\,\E_\zs{y}\,\prod^{n-1}_\zs{j=1}\,
g_\zs{j+1}(\Phi_\zs{j})\
\mu(\d y)
=
\int_\zs{\cX}\,
\E_\zs{x}\,G_\zs{n}(\Phi_\zs{1},\ldots,\Phi_\zs{n})
\,\lambda(\d x)\,.
$$
Hence, the Proposition~\ref{Pr.subsec:AMc.5}.
\endproof

\begin{proposition}\label{Pr.subsec:AMc.6}
Assume that the splitting chain
$(\check{\Phi}_\zs{n})_\zs{n\ge 1}$
has an invariant probability measure $\check{\pi}$. Then,
the chain $(\Phi_\zs{n})_\zs{n\ge 1}$ has the invariant probability
measure $\pi$ on $\cB(\cX)$ which is given as
\begin{equation}\label{subsec:AMc.14}
\pi(\Gamma)=\check{\pi}(\Gamma_\zs{0})+
\check{\pi}(\Gamma_\zs{1})\,.
\end{equation}
Moreover, $\check{\pi}=\pi^{*}$.
\end{proposition}
\proof First we check directly that $\check{\pi}=\pi^{*}$. Moreover,
for any $\Gamma\in\cB(\cX)$
\begin{align*}
\pi(\Gamma)&=\check{\pi}(\Gamma_\zs{0}\cup \Gamma_\zs{1})
=\int_\zs{\check{\cX}}\,
\check{\pi}(\d \check{z})\,\check{\P}^{\vartheta}(\check{z},\Gamma_\zs{0}\cup \Gamma_\zs{1})\\[2mm]
&
=\int_\zs{\check{\cX}}\,
\pi^{*}(\d \check{z})\,\check{\P}^{\vartheta}(\check{z},\Gamma_\zs{0}\cup \Gamma_\zs{1})
\,.
\end{align*}
Therefore, applying here Proposition~\ref{Pr.subsec:AMc.5}
we obtain that
$$
\pi(\Gamma)=\int_\zs{\cX}\,\P(z,\Gamma)\,\pi(\d z)\,,
$$
i.e. $\pi$ is the invariant measure for the chain
$(\Phi_\zs{n})_\zs{n\ge 1}$. Hence
Proposition~\ref{Pr.subsec:AMc.6}.
\endproof

\medskip

\subsection{Moment inequality for the process (\ref{sec:In.3})}

\begin{proposition}\label{Pr.subsec:AMc.7}
Let $(y_t)_\zs{t\ge 0}$ be a solution of the equation \eqref{sec:In.3}. Then, for any $m\ge 1$
and for any $x\in\bbr$,
\begin{equation}\label{subsec:AMc.15}
\sup_\zs{t\ge 0}\,
\sup_\zs{\vartheta\in\Theta}\,
\E^{\vartheta}_\zs{x}\,(y_\zs{t})^{2m}\,\le (2m-1)!!\,(x^{2}+\M_\zs{*}/\beta)^{m}
\,,
\end{equation}
where $\M_\zs{*}=
(M+\beta\x_\zs{*})^2/\beta+2\sigma^{2}_\zs{max}$.
\end{proposition}
\proof To obtain this inequality we make use of the method proposed
in (\cite{KaPe}, p.20) for linear stochastic equation. First of
all note that thanks to Theorem 4.7 from \cite{LpSh}, for any
$T>0$, there exists some $\epsilon>0$ such that for each
$\vartheta\in\Theta$ and $x\in\bbr$
\begin{equation}\label{subsec:AMc.15-1}
\sup_\zs{0\le t\le T}\, \E^{\vartheta}_\zs{x}\,
e^{\epsilon y^{2}_\zs{t}}\, <\infty \,.
\end{equation}
Applying the Ito formula to $y_\zs{t}^{2m}$ and denoting
$$
\B_\zs{\vartheta}(y)=2yS(y)+\sigma^{2}(y)+\beta y^{2}
\,,
$$
yield
\begin{align*}
\d y_\zs{t}^{2m}&=-m\beta y_\zs{t}^{2m}\d t+
m y_\zs{t}^{2(m-1)}\left(\B_\zs{\vartheta}(y_\zs{t})+2(m-1)\sigma^{2}(y_\zs{t}) \right)\d t
\\[2mm]
&+2m y_\zs{t}^{2m-1} \sigma(y_\zs{t}) \d W_\zs{t}\,.
\end{align*}
Therefore, taking into account that $y_\zs{0}=x$ we can represent the last equation 
in the following integral form
\begin{align}\nonumber
y^{2m}_\zs{t}&=e^{-m\beta t}\,x^{2m}+m
\int^{t}_\zs{0}\,e^{-m\beta (t-s)}\,
y_\zs{s}^{2(m-1)}\left(\B_\zs{\vartheta}(y_\zs{s})+2(m-1)\sigma^{2}(y_\zs{s}) \right)\d s
\\[2mm]\label{subsec:AMc.15-2}
&+2m\int^{t}_\zs{0}\,e^{-m\beta(t-s)}\,y_\zs{s}^{2m-1} \sigma(y_\zs{s}) \d W_\zs{s}\,.
\end{align}
One can check directly that
$$
\sup_\zs{\vartheta\in\Theta}\,
\sup_\zs{y\in\bbr}\,
\,|\B_\zs{\vartheta}(y)|\,\le 
\frac{
(M+\beta\x_\zs{*})^2}{\beta}
+\sigma^{2}_\zs{max}\,.
$$
Moreover, the property \eqref{subsec:AMc.15-1} yields that, for any $m\ge 1$,
$$
\E^{\vartheta}
\int^{t}_\zs{0}\,
e^{-m\beta(t-s)}\,
y_\zs{s}^{2m-1} \sigma(y_\zs{s}) \d W_\zs{s}=0\,.
$$
Therefore, setting $z_\zs{t}(m)=\E^{\vartheta}_\zs{x} y^{2m}_\zs{t}$,
 we obtain 
$$
z_\zs{t}(m)\le x^{2m}+m(2m-1)\M_\zs{*}
\,\int^{t}_\zs{0}\,e^{-m\beta (t-s)}\,z_\zs{s}(m-1)\,\d s\,.
$$
The induction implies directly the bound \eqref{subsec:AMc.15}.
Hence Proposition~\ref{Pr.subsec:AMc.7}.
\endproof

\begin{proposition}\label{Pr.subsec:AMc.8}
Let $(y_t)_\zs{t\ge 0}$ be a solution of the equation \eqref{sec:In.3}.
Then, for any 
$K>\sqrt{\M_\zs{1}}$,
\begin{equation}\label{subsec:AMc.16}
\sup_\zs{|x|\le K}\,
\sup_\zs{\vartheta\in\Theta}\,
\P^{\vartheta}_\zs{x}\,
\left(|y_\zs{1}|\ge K\right)\,
\le
\frac{4\sigma^{2}_\zs{max}
\left(K^{2}+\M_\zs{2}\right)}{\beta(1-e^{-\beta})\left(K^{2}-\M_\zs{1}\right)^{2}}
\,,
\end{equation}
where $\M_\zs{2}$ and 
$\M_\zs{1}$ are given in \eqref{sec:De.6-3}.
\end{proposition}
\proof First, putting in \eqref{subsec:AMc.15-2} $m=1$,
we obtain
$$
\sup_\zs{t\ge 0}
\E^{\vartheta}_\zs{x}\,y^{2}_\zs{t}\,
\le\,x^{2}+\M_\zs{2}
$$
and
$$
\P^{\vartheta}_\zs{x}\,
\left(y^{2}_\zs{1}\ge K^{2}\right)\le
\P\left( 2\zeta
\ge \left(K^{2}-\M_\zs{1}\right)(1-e^{-\beta})
\right)\,,
$$
where $\zeta=\int^{1}_\zs{0}\,e^{-\beta(1-s)}\,y_\zs{s} \sigma(y_\zs{s}) \d W_\zs{s}$.
Now, taking into account that for $|x|\le K$
$$
\E^{\vartheta}_\zs{x}\,\zeta^{2}=
\int^{1}_\zs{0}\,e^{-2\beta(1-s)}\,\E^{\vartheta}_\zs{x}\,y^{2}_\zs{s} \sigma^{2}(y_\zs{s}) \d s
\le \sigma^{2}_\zs{max}\left(
K^{2}+\M_\zs{2}
\right)
\frac{1-e^{-\beta}}{\beta}\,.
$$
The Chebyshev inequality implies now the bound \eqref{subsec:AMc.16}. Hence
Proposition~\ref{Pr.subsec:AMc.8}. \endproof

\end{document}